\documentclass{gtart_h}

\def\ifplaintex{\expandafter\ifx\csname documentclass\endcsname\relax}

\def\gtp{{\mathsurround=0pt\it $\cal G\mskip-2mu$eometry \&\ 
$\cal T\!\!$opology $\cal P\!$ublications}}  

\def\recd{{\small Received:\qua\receiveddate\ifx\reviseddate\relax
\else\qquad Revised:\qua\reviseddate\fi\par}} 

\def\lognumber#1{\def\thelognumber{#1}}
\def\volumenumber#1{\def\thevolumenumber{#1}}
\def\volumeyear#1{\def\thevolumeyear{#1}}
\def\papernumber#1{\def\thepapernumber{#1}}
\def\pagenumbers#1#2{\def\startpage{#1}\def\finishpage{#2}}
\def\published#1{\def\publishdate{#1}}

\def\received#1{\def\receiveddate{#1}}
\def\revised#1{\def\reviseddate{#1}}
\def\accepted#1{\def\accepteddate{#1}}

\def\asciiaddress#1{\def\theasciiaddress{#1}}

\long\def\asciiabstract#1{\long\def\theasciiabstract{#1}}

\let\\\par\let\thelognumber\relax\let\thevolumenumber\relax
\let\thepapernumber\relax\let\thevolumeyear\relax\let\startpage\relax
\let\finishpage\relax\let\publishdate\relax\let\receiveddate\relax
\let\reviseddate\relax\let\accepteddate\relax\let\theasciititle\relax
\let\theasciiauthors\relax\let\theasciiaddress\relax
\let\theasciiabstract\relax

\let\theasciiemail\relax

\ifplaintex
\font\logobig=cmssbx10 scaled 3836
\font\logomed=cmssbx10 scaled 2557
\else
\font\logobig=cmssbx10 scaled 4200
\font\logomed=cmssbx10 scaled 2800
\fi

\long\def\makeagttitle{   
\count0=\startpage
\agt\hfill      
\hbox to 45truept{\vbox to 0pt{\vglue -13truept{\logomed A\kern -.37em{\logobig 
T}\kern -.38em G}\vss}\hss}
\break
{\small Volume \thevolumenumber\ (\thevolumeyear)
\startpage--\finishpage\nl
Published: \publishdate}

\vglue .25truein

{\parskip=0pt\leftskip 0pt plus
1fil\def\\{\par\smallskip}{\Large\bf\thetitle}\par\medskip} \vglue
0.05truein

{\parskip=0pt\leftskip 0pt plus 1fil\def\\{\par}{\sc\theauthors}
\par\medskip}%
 
\vglue 0.03truein 

{\small\leftskip 25truept\rightskip 25truept{\bf Abstract}\stdspace\theabstract

{\bf AMS Classification}\stdspace\theprimaryclass
\ifx\thesecondaryclass\relax\else; \thesecondaryclass\fi\par
{\bf Keywords}\stdspace \thekeywords\par}\vglue 7truept

}   

\ifplaintex
\font\phead=cmsl9 scaled 950
\font\pnum=cmbx10 scaled 913
\font\pfoot=cmsl9 scaled 950
\headline{\vbox to 0pt{\vskip -4.5mm\line{\small\phead\ifnum
\count0=\startpage ISSN 1472-2739 (on-line) 1472-2747 (printed)
\hfill {\pnum\folio}\else\ifodd\count0\def\\{ }%
\ifx\theshorttitle\relax\thetitle\else\theshorttitle\fi\hfill{\pnum\folio}
\else\def\\{ and }{\pnum\folio}\hfill\ifx\theshortauthors\relax\theauthors
\else\theshortauthors\fi\fi\fi}\vss}}
\footline{\vbox to 0pt{\vglue 0mm\line{\small\pfoot\ifnum\count0=\startpage
\copyright\ \gtp\hfill\else
\agt, Volume \thevolumenumber\ (\thevolumeyear)\hfill\fi}\vss}}
\else
\font=cmsl9 scaled 1050
\font\lnum=cmbx10 
\font=cmsl9 scaled 1050
\makeatletter
\def\@oddhead{{\small\lhead\ifnum\count0=\startpage ISSN 1472-2739 
(on-line) 1472-2747 (printed)\hfill {\lnum\number\count0}\else\ifodd\count0
\def\\{ }\ifx\theshorttitle\relax \thetitle \else\theshorttitle\fi\hfill
{\lnum\number\count0}\else\def\\{ and }{\lnum\number\count0}
\hfill\ifx\theshortauthors\relax 
\theauthors\else\theshortauthors\fi\fi\fi}}\def\@evenhead{\@oddhead}
\def\@oddfoot{\small\lfoot\ifnum\count0=\startpage\copyright\ \gtp\hfill\else
\agt, Volume \thevolumenumber\ (\thevolumeyear)\hfill\fi}
\def\@evenfoot{\@oddfoot}
\makeatother
\fi
\let\maketitlepage\makeagttitle
\let\makeshorttitle\maketitlepage
\let\maketitle\maketitlepage

\newwrite\gtoutfile
\long\gdef\makeheadfile{  
{\def\\{, }\def\s{ }
\immediate\openout\gtoutfile head.xxx
\immediate\write\gtoutfile{Proxy-for: \ifx\theasciiauthors\relax
\theauthors\else\theasciiauthors\fi\s<\ifx\theasciiemail\relax\theemail\else\theasciiemail\fi>}
\immediate\write\gtoutfile{\noexpand\\}
\immediate\write\gtoutfile{Authors: \ifx\theasciiauthors\relax
\theauthors\else\theasciiauthors\fi}
{\def\\{ }\immediate\write\gtoutfile{Title: \ifx\theasciititle\relax
\thetitle\else\theasciititle\fi}}
\immediate\write\gtoutfile{Subj-class: GT or SG, GR etc}
\immediate\write\gtoutfile{MSC-class: \theprimaryclass\ifx\thesecondaryclass\relax\else, \thesecondaryclass\fi}
\immediate\write\gtoutfile{Journal-ref: Algebr. Geom. Topol. \thevolumenumber\s
(\thevolumeyear) \startpage-\finishpage}
\immediate\write\gtoutfile{Comments: Published by Algebraic and
Geometric Topology at}
\immediate\write\gtoutfile{\s\s\s  http://www.maths.warwick.ac.uk/agt/AGTVol\thevolumenumber/agt-\thevolumenumber-\thepapernumber.abs.html}
\immediate\write\gtoutfile{\noexpand\\}
\immediate\write\gtoutfile{}
\ifx\theasciiabstract\relax
\immediate\write\gtoutfile{\theabstract}\else
\immediate\write\gtoutfile{\theasciiabstract}\fi
\immediate\write\gtoutfile{}
\immediate\write\gtoutfile{\noexpand\\}
\immediate\write\gtoutfile{}
\immediate\closeout\gtoutfile}}  

\def\maketitlepage{\makeagttitle\makeheadfile}
\let\makeshorttitle\maketitlepage
\let\maketitle\maketitlepage

\lognumber{53}
\volumenumber{4}
\volumeyear{2004}
\papernumber{53}
\pagenumbers{1211}{1251}
\received{24 January 2004} 
\revised{18 November 2004}
\accepted{8 December 2004}
\published{21 December 2004}

\usepackage{graphicx,psfrag,amssymb,amsmath,amscd,xypic}

\theoremstyle{plain}
\newtheorem{theorem}{Theorem}
\newtheorem*{theorem*}{Theorem}
\newtheorem{proposition}{Proposition}[section]
\newtheorem{lemma}{Lemma}[section]

\theoremstyle{definition}
\newtheorem{definition}{Definition}
\newtheorem*{conjecture}{Conjecture}

\theoremstyle{remark}
\newtheorem*{remark}{Remark}

\begin{document}
\title{An invariant of link cobordisms\\from Khovanov homology}
\author{Magnus Jacobsson}
\address{Istituto Nazionale di Alta Matematica
(INdAM), Citt\`a Universitaria\\P.le Aldo Moro 5, 00185 Roma,
Italy} 
\asciiaddress{Istituto Nazionale di Alta Matematica
(INdAM), Citta Universitaria\\P.le Aldo Moro 5, 00185 Roma,
Italy} 
\email{jacobsso@mat.uniroma1.it}

\begin{abstract} In \cite{K}, Mikhail Khovanov constructed a homology theory 
for oriented links, whose graded Euler characteristic is the 
Jones polynomial. He also explained how every link
cobordism between two links induces a homomorphism between their
homology groups, and he conjectured the invariance 
(up to sign) of this homomorphism under ambient isotopy of the 
link cobordism. In this paper we prove this conjecture, after 
having made a necessary improvement on its statement. 
We also introduce polynomial Lefschetz numbers of cobordisms 
from a link to itself such that the Lefschetz polynomial 
of the trivial cobordism is the Jones polynomial. These
polynomials can be computed on the chain level.
\end{abstract}
\asciiabstract{%
In [Duke Math. J. 101 (1999) 359-426], Mikhail Khovanov constructed a
homology theory for oriented links, whose graded Euler characteristic
is the Jones polynomial.  He also explained how every link cobordism
between two links induces a homomorphism between their homology
groups, and he conjectured the invariance (up to sign) of this
homomorphism under ambient isotopy of the link cobordism. In this
paper we prove this conjecture, after having made a necessary
improvement on its statement.  We also introduce polynomial Lefschetz
numbers of cobordisms from a link to itself such that the Lefschetz
polynomial of the trivial cobordism is the Jones polynomial. These
polynomials can be computed on the chain level.}

\primaryclass{57Q45}
\secondaryclass{57M25}
\keywords{Khovanov homology, link cobordism, Jones polynomial}
\makeshorttitle


\section{Introduction}
\label{i}
\noindent
In \cite{K}, M Khovanov associated to any diagram $D$ of an oriented 
link a chain complex $C(D)$ of abelian groups, whose Euler 
characteristic is the Jones polynomial \cite{J}. He proved that for 
any two diagrams of the same link the corresponding 
complexes are chain equivalent. Hence, the 
homology groups $\mathcal{H}(D)$ are link invariants up to isomorphism.
For a definition of the chain complex, see Definitions \ref{Dcijl}
through \ref{Ddiff}, Section \ref{cg} below.  See also Bar-Natan \cite{BN}
for a treatment of Khovanov homology.

One of the motivations of Khovanov's work was the hope of finding a
lift of the Penrose--Kauffman spin networks calculus to a calculus 
of surfaces in the 4--sphere. To finish this program he suggested 
the following TQFT construction of an invariant of link cobordisms. 

Any link cobordism can be described as a one-parameter family 
$D_t, t \in [0,1]$ of planar diagrams, called a {\em movie}. 
The $D_t$ are link diagrams, except at finitely many $t$--values 
where the topology changes: the diagram undergoes a {\em local move}, 
which is either a Reidemeister move or a Morse modification. 
Away from these values the diagram experiences a planar isotopy as $t$ 
varies. Khovanov explained how local moves induce chain 
maps between complexes, hence homomorphisms between 
homology groups. The same is true for planar isotopies. 
Hence, the composition of these chain maps defines a homomorphism 
between the homology groups of the diagrams of the boundary links.

Khovanov conjectured that, up to multiplication by $-1$, this homomorphism 
is invariant under ambient isotopy of the link cobordism. 
(For the exact formulation see Section \ref{ce}.)

In this paper, we show that this conjecture is not properly 
stated, by giving simple counterexamples (Theorem \ref{818}, 
Section \ref{cekc}). In fact, there are very simple 
examples for multi-component links, but we also 
give an example in the case of knots. We then show that this 
can be remedied by considering only ambient isotopies which 
leave the links in the boundary setwise fixed. If the 
conjecture is modified in this way, it is indeed true.
This is the main result of the paper.

\begin{theorem*}[Theorem \ref{mkc}, Section \ref{rtc}]
For oriented links $L_0$ and $L_1$, presented by diagrams $D_0$ and $D_1$, 
an oriented link cobordism $\Sigma$ from $L_0$ to $L_1$, defines
a homomorphism $\mathcal{H}(D_0) \rightarrow \mathcal{H}(D_1)$, invariant 
up to multiplication by -1 under ambient isotopy of $\Sigma$ leaving
$\partial \Sigma$ setwise fixed. Moreover, this invariant is non-trivial.
\end{theorem*} 

The proof of Khovanov's conjecture implies the existence 
of a family of derived invariants of link cobordisms with the same source 
and target, which are analogous to the classical Lefschetz numbers
of endomorphisms of manifolds. We call these {\em Lefschetz polynomials 
of link endocobordisms} and like their classical analogues 
they are computable on the chain level. The Jones polynomial 
appears as the Lefschetz polynomial of the identity cobordism 
(Section \ref{lple}).

A knotted closed surface is a link cobordism between empty 
links. The grading properties of the theory force the 
invariant to be zero on such a surface unless its 
Euler characteristic is zero (see Section \ref{def}). 

To get a theory which does not immediately rule out knotted 
spheres, Khovanov actually stated his conjecture in a more general 
setting, defining the chain complex as a module over the polynomial ring 
$\mathbb{Z}[c]$, with a differential also depending on $c$. The 
integer theory is obtained by setting $c=0$. Surprisingly, 
it turns out that the conjectured invariant with 
polynomial coefficients does not exist. A proof of this fact 
is given in \cite{Ja2}.

The following is an outline of the paper. Section \ref{kh} contains 
an elementary description of Khovanov homology, suggested by Oleg 
Viro \cite{V}. In Section \ref{lc} we introduce link cobordisms 
and their diagrams, and explain how they induce maps 
between homology groups. Section \ref{ce} contains the 
counterexamples previously mentioned, and Section \ref{pmkc} 
contains the main result of the paper,
the proof of Khovanov's conjecture. Finally, the last section 
introduces the Lefschetz polynomials.


\section{Khovanov homology}
\label{kh}
In this section we give an elementary description of Khovanov's 
homology theory. 


\subsection{Chain groups}
\label{cg}

\begin{figure}[ht!]\small
\begin{center}
\psfrag{positive}{positive}
\psfrag{negative}{negative}
\includegraphics[width = 4 cm, height = 4 cm]{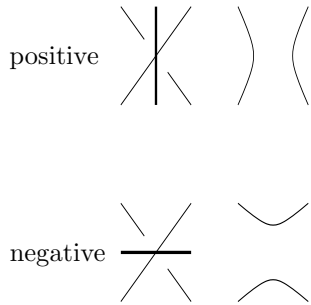}
\caption{Resolution of a state according to markers}
\label{resolution}
\end{center}
\end{figure}

Let $D$ be an oriented link diagram and $a$ one of its crossings. 
A {\em marker at a} is a choice of a pair of vertical
angles at $a$. In pictures this choice is indicated by a short bar
connecting the chosen angles (see Figure \ref{resolution}).

If the bar is in the region swept out by the overcrossing line 
as it is turned counterclockwise toward the undercrossing line, 
the marker is called positive. Otherwise, it is negative.

\begin{figure}[ht!]\small
\begin{center}
\psfrag{+}{$+$}
\psfrag{-}{$-$}
\includegraphics[width = 6 cm, height = 3 cm]{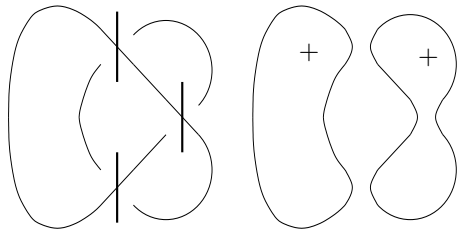}
\caption{A state of a diagram of the unknot: here $i(S) = 0$, $j(S) = -3$. 
Recall that the writhe of a knot (but not of a link) is 
independent of the choice of orientation.}
\label{atrefoil}
\end{center}
\end{figure}

To a distribution $s$ of markers over the crossings of $D$ corresponds
an embedded collection $C_1,...,C_{r_s}$ of circles in the plane, called 
the {\em resolution of $D$ according to $s$}. It is obtained by
replacing a small neighbourhood of each crossing point by a 
pair of parallel arcs in the regions not specified by the marker 
(see Figure \ref{resolution}). A {\em signed resolution} is a resolution 
together with a choice of $+/-$--sign for each circle $C_i$. 
A {\em state of D} is a distribution of markers at the crossing points,
together with a signed resolution of it. The signed resolution of a
state $S$ will be denoted by res($S$) (see Figure \ref{atrefoil}).

Let $S$ be a state. Denote by $\sigma(S)$ the sum of the signs
of the markers in $S$ and by $\tau(S)$ the sum of the signs in its
signed resolution res($S$). Furthermore, let $w(D)$ denote the writhe number
of the diagram $D$. The following functions will be the grading
parameters of Khovanov's chain complex.
\begin{displaymath}
\begin{split}
i(S) &= \frac{w(D) - \sigma(S)}{2} \\
j(S) &= - \frac{\sigma(S) + 2\tau(S) - 3w(D)}{2} 
\end{split}
\end{displaymath}

\begin{remark}
The definition of a state of a diagram is a refinement of L Kauffman's 
definition, which was used in \cite{Ka} to construct the Jones polynomial 
$V_L$ of a link $L$. The refinement consists of the signs 
associated to the components of the resolution. In \cite{Ka} a
state is only a distribution of markers together with the corresponding 
resolution. 
To understand this refinement, consider the equalities below. 
\begin{equation*}
\begin{split}
V_L &= \sum_{\textrm{Kauffman states s of D}} A^{\sigma(s)} (-A)^{-3w(D)} (-A^2 -
A^{-2})^{r_s} \\
&= \sum_{\textrm{``refined'' states S of D}} (-1)^{w(D)} A^{\sigma(S) - 3w(D)}
(-A^{-2})^{- \tau(S)} \\
&= \sum_{\textrm{``refined'' states S of D}} (-1)^{\frac{w(D)-\sigma(S)}{2}}
q^{-\frac{\sigma(S)+2\tau(S)-3w(D)}{2}} \\ 
&=\sum_{\textrm{``refined'' states S of D}} (-1)^{i(S)} q^{j(S)} \\
\end{split}
\end{equation*}
The first equality is Kauffman's definition. Each term in the sum 
corresponds to an ``unrefined'' state $s$. Each term is polynomial and
contains a factor $(-A^2 - A^{-2})^{r_s}$. If we associate the term 
$-A²$ with positive circles and the term $-A^{-2}$ with negative
circles, then the refinement identifies each monomial in the binomial 
expansion of this polynomial term with a state in our sense. 
Letting $q=-A^{-2}$, we get $V_L(q)$ as a sum over the refined states.
\end{remark}

\begin{definition}
\label{Dcijl}
Let $L$ be a subset of the set $I$ of crossings of $D$.
Let $C_{L}^{i,j}(D)$ be the free abelian group on the set of states $S$
which have $i(S) = i$, $j(S) = j$ and $L$ as the set of crossings 
with negative markers.
\end{definition} 

\begin{remark} Observe that the cardinality $|L| = n_i$ of $L$
is uniquely determined by $i$. 
\end{remark}

For any finite set $S$, let $FS$ be the free abelian group generated by $S$. 
For bijections $f,g:\{1,...|S|\} \rightarrow S$, let $p(f,g) \in
\{0,1\}$ be the parity of the permutation $f^{-1}g$ of
$\{1,...|S|\}$. Let $Enum(S)$ be the set of all such bijections 
$f,g$. 

\begin{definition}
For $S$ as above, we define
\begin{displaymath}
E(S) = FEnum(S)/((-1)^{p(f,g)} f - g).
\end{displaymath}
\end{definition}

Khovanov's chain groups are defined as follows.

\begin{definition}
The $(i,j)$-th chain group of the chain complex is
\begin{displaymath}
C^{i,j}(D) = \bigoplus_{L \subset I, \mid L \mid = n_i} C^{i,j}_{L}(D) 
\otimes E(L).
\end{displaymath}
\end{definition}


\subsection{The differential}
\label{d}

To define the differential, note that replacing a positive marker in a
state $S$ with a negative one changes the resolution in one of two possible
ways; the number of circles either increases or decreases by one. 

\begin{figure}[ht!]\small
\begin{center}
\psfrag{p}{$p$}
\psfrag{q}{$q$}
\psfrag{p:q}{$p{:}q$}
\psfrag{q:p}{$q{:}p$}
\psfrag{0}{0}
\includegraphics[width = 7 cm, height = 10 cm]{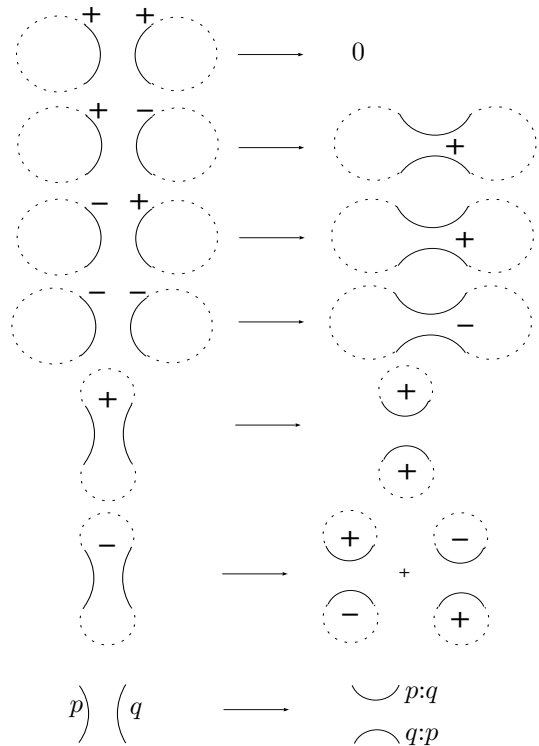}
\caption{A Frobenius calculus of signed circles}
\label{frobenius}
\end{center}
\end{figure}

\begin{definition}
\label{Ddiff}
Let $S$ belong to $C^{i,j}_{L}(D)$. The differential of $S \otimes [x]$
is the sum
\begin{displaymath}
d(S \otimes [x]) = \sum_{T} T \otimes [xa], 
\end{displaymath}
where the $T$:s run over all states in $C^{i+1,j}(D)$ which satisfy 
the restrictions in the itemized list below, determined by $S$. 
Here, $x$ is an ordered sequence of crossings and $a = a(T)$ 
a special crossing explained below. Square brackets around 
a sequence of crossings denote its equivalence class in $E(L)$.

\begin{itemize}
\item $T$ has the same markers as $S$ except at one crossing point
$a$. At this point $T$ has a negative marker, whereas $S$ 
has a positive one. 
\item This restriction on markers implies that the resolutions of 
$T$ and $S$ coincide outside a small disc neighbourhood of $a$. 
We require the signed resolution res($T$) to coincide with res($S$) 
on the circles which do not intersect this disc. 
The disc is intersected by either one or two components of 
res($S$). In the former case res($T$) intersects the 
disc with two of its components; in the latter case, one.
We require the total sum of signs of the circles of $T$ to be
greater by one than that of $S$ (see Figure \ref{frobenius} and the
remark below).
\end{itemize}
\end{definition}
It is straightforward to check that the map $d$ thus defined is
actually a differential (i.e.\ that $d^2 = 0$) and that the bidegree
of $d$ is $(1,0)$.

\begin{remark}
When both arcs passing the crossing $a$ belong to the 
{\em same negative} component of res($S$), 
there are two ways to distribute signs to res($T$) in accordance 
with the restrictions. This explains why (and in what sense) the sixth 
row of Figure \ref{frobenius} is a sum.
When there are {\em two positive components} of $S$ passing $a$, the
resolution of $T$ cannot be supplied with signs consistent 
with the restrictions. Hence the zero in the first row.

When signs $+$, $-$ are used in figures, two arcs 
which belong to the same component are both labelled with one 
common, single sign, as in the first six rows of Figure
\ref{frobenius}. 

However, we will also use the mnemonic notation 
on the last row of Figure \ref{frobenius}. It summarizes the
other rows in the following way. The labels $p$, $q$, $p{:}q$, $q{:}p$ 
represent the signs. Thus, 
in this notation there are always four labels on four arcs,
regardless of whether the arcs are connected or not outside 
the figure. Furthermore, the labels $p{:}q$, $q{:}p$ on arcs 
in a figure may mean that the state appears with coefficient zero 
(if the arcs labelled $p$ and $q$ are in different and positive 
components, as in the first row), or that the sum of two states 
appears, with opposite signs on the two arcs (if the arcs labelled 
$p$ and $q$ are in the same, negative, component).
\end{remark}
\begin{remark}
One may regard the signs $+$ and $-$ as linear generators of a commutative 
Frobenius algebra $A$, with multiplication and comultiplication defined by the
calculus in the previous remark.
It is a theorem (see e.g.\ \cite{A}) 
that commutative Frobenius algebras are in one-to-one correspondence 
with $(1+1)$--dimensional topological quantum field theories. 

Indeed, this is Khovanov's starting point in \cite{K}. 
To each collection of $k$ circles in the plane
his TQFT associates the tensor product of $k$ copies of $A$. 
Any compact oriented surface with boundary the disjoint union of two
such collections (with cardinalities $k$,$n$), induces a map 
$A^{\otimes k} \rightarrow A^{\otimes n}$, via its decomposition
into discs and ``pairs of pants'', and the following descriptions.

To any pair of pants directed from the feet upwards corresponds 
the product in the algebra, $A\otimes A \xrightarrow{m} A$.
To a pair of pants directed from the waist downwards corresponds
the coproduct $A \xrightarrow{\Delta} A \otimes A$. 
A disc induces either a counit $A \xrightarrow{e} \mathbb{Z}$ or a unit
$\mathbb{Z} \xrightarrow{i} A$, depending on whether it is regarded
as directed from the empty set to the circle, or conversely.
It follows from the axioms in the Frobenius algebra that the induced map 
is independent of the decomposition.

The values of multiplication and comultiplication 
on generators are given by the rules in Figure
\ref{frobenius} in the following sense:
\begin{displaymath}
\begin{split}
m(+ \otimes +) &= 0 \\
m(+ \otimes -) &= m(- \otimes +) = + \\
m(- \otimes -) &= - \\
\Delta(+) &= + \otimes +\\
\Delta(-) &= (+ \otimes -) + (- \otimes +) 
\end{split}
\end{displaymath}
The maps for discs are given by:
\begin{displaymath}
\begin{split}
i(1) &= - \\
e(+) &= 1 \\
e(-) &= 0 \\
\end{split}
\end{displaymath} 
\end{remark}


\subsection{Khovanov homology}
As mentioned above, the differential $d$ has bidegree $(1,0)$. Hence, 
for each $j$, $C^{i,j}(D)$ is an ordinary chain complex 
graded by $i$. The homology $\mathcal{H}^{i,j}(D)$ of this complex 
has an Euler characteristic which can be computed as the 
alternating sum of the ranks of the chain groups. Summing 
over $j$ with coefficients $q^j$ we obtain
\begin{displaymath}
\sum_{i,j}(-1)^iq^j\operatorname{rk}\mathcal{H}^{i,j}(D) = 
\sum_{i,j}(-1)^iq^j\operatorname{rk}C^{i,j}(D). 
\end{displaymath}
This is known as {\em the graded Euler characteristic} 
of the complex $C^{i,j}(D)$ (or of $\mathcal{H}^{i,j}(D)$). 
From the definition of $i$ and $j$, it follows that the Jones polynomial 
is the graded Euler characteristic of Khovanov's chain complex.


\section{Link cobordisms} 
\label{lc}


\subsection{Link cobordisms and their diagrams}
\label{lctd}
\begin{definition}
Let $I = [0,1]$. An {\em (oriented) link cobordism} between two links 
$L_0 \subset \mathbb{R}^3 = \mathbb{R}³ \times \{0\}$ and $L_1 \subset
\mathbb{R}³ = \mathbb{R}^3 \times \{1\}$ is a smooth, compact, oriented 
surface neatly embedded in $\mathbb{R}^3 \times I$ with boundary 
$\partial \Sigma = $$\Sigma \cap (\mathbb{R}^3 \times \partial{I})$$ = L_0 
\cup L_1$. Close to the boundary, $\Sigma$ is orthogonal to 
$\mathbb{R}^3 \times \partial{I}$.
\end{definition}

\begin{remark}We will refer to $L_0$ and $L_1$ as the {\em source} 
and {\em target} of the cobordism $\Sigma$. 
\end{remark}

The orientation of $\Sigma$ induces an orientation on all the  
intersections $L_t$ of $\Sigma$ with the constant time hyperplanes 
$\mathbb{R}^3 \times \{t\}$. A tangent vector $v$ to a link component 
is in the positive direction if $(v,w)$ gives the orientation of $\Sigma$ 
whenever $w$ is a vector tangent to $\Sigma$ in the direction of 
increasing time. 

\begin{definition}
A link cobordism $\Sigma$ is {\em generic} if time $t$ is a Morse function 
on $\Sigma$ with distinct critical values. 
\end{definition}

A generic link cobordism intersects hyperplanes of constant 
$t \in [0,1]$ in embedded links except for a finite set 
of values $t$, for which the intersection $L_t$ either has a single 
transversal double point or is the disjoint union of a link and 
an isolated point.

A generic link cobordism can be represented by a three-dimensional 
{\em surface diagram}, directly analogous to the two-dimensional 
diagrams of classical links (or tangles). Such a diagram is the 
image of the link cobordism under a projection to $\mathbb{R}^2 
\times [0,1]$ which preserves the $t$--variable. The projection is 
required to be generic in the sense that the only singular points in 
the interior of the surface diagram are double points, Whitney umbrella 
points and triple points. At a double point, the diagram 
looks like the transversal intersection of two planes. Whitney
umbrellas and triple points occur as the (isolated) boundary 
points of the double point set in the interior of $\mathbb{R}^2 
\times [0,1]$.  

A full set of Reidemeister-type moves for surface diagrams was
given by D Roseman \cite{R} and are called {\em Roseman moves}.
These moves (plus ambient isotopy of the diagram) are sufficient 
to transform any two surface diagrams of ambient isotopic 
surfaces into each other. We will use another way of presenting 
a link cobordism. 

\begin{definition} A {\em movie} of a generic link 
cobordism with a given surface diagram as above, is the 
intersection of the diagram 
with planes $\mathbb{R}^2 \times \{t\}$, regarded
as a function of time $t$. The intersection for a fixed $t$ is called a
{\em still}. The (finitely many) $t$ for which the intersection 
is not a link diagram are called {\em critical levels}. 
\end{definition}

Observe that the restriction of the movie to a small interval of 
time around a critical level shows a link diagram which undergoes 
either a Reidemeister move or an oriented Morse modification. 
The Reidemeister moves occur at (interior) levels where the 
double point set has a boundary point or a local maximum 
or minimum. The Morse modifications occur at smooth points of 
the surface diagram which are minima, maxima or saddle points
 of the time function. Between two critical levels the diagram 
undergoes a planar isotopy.

\begin{definition} Reidemeister moves and Morse modifications will be
called {\em local moves}. Each such move is localized to a small
disc, called a {\em changing disc}.
\end{definition}

When a movie is studied, not all the stills of it are 
considered, but only as many as are necessary for a 
picture of the surface; one still between each pair 
of consecutive critical points and stills of the source 
and target diagrams (and sometimes a few additional ones, 
for clarity). 

\subsection{Movie moves}
\label{mm}
\begin{figure}[ht!]\small
\begin{center}
\psfrag{1}{1}
\psfrag{2}{2}
\psfrag{3}{3}
\psfrag{4}{4}
\psfrag{5}{5}
\includegraphics[width = 13 cm, height = 7 cm]{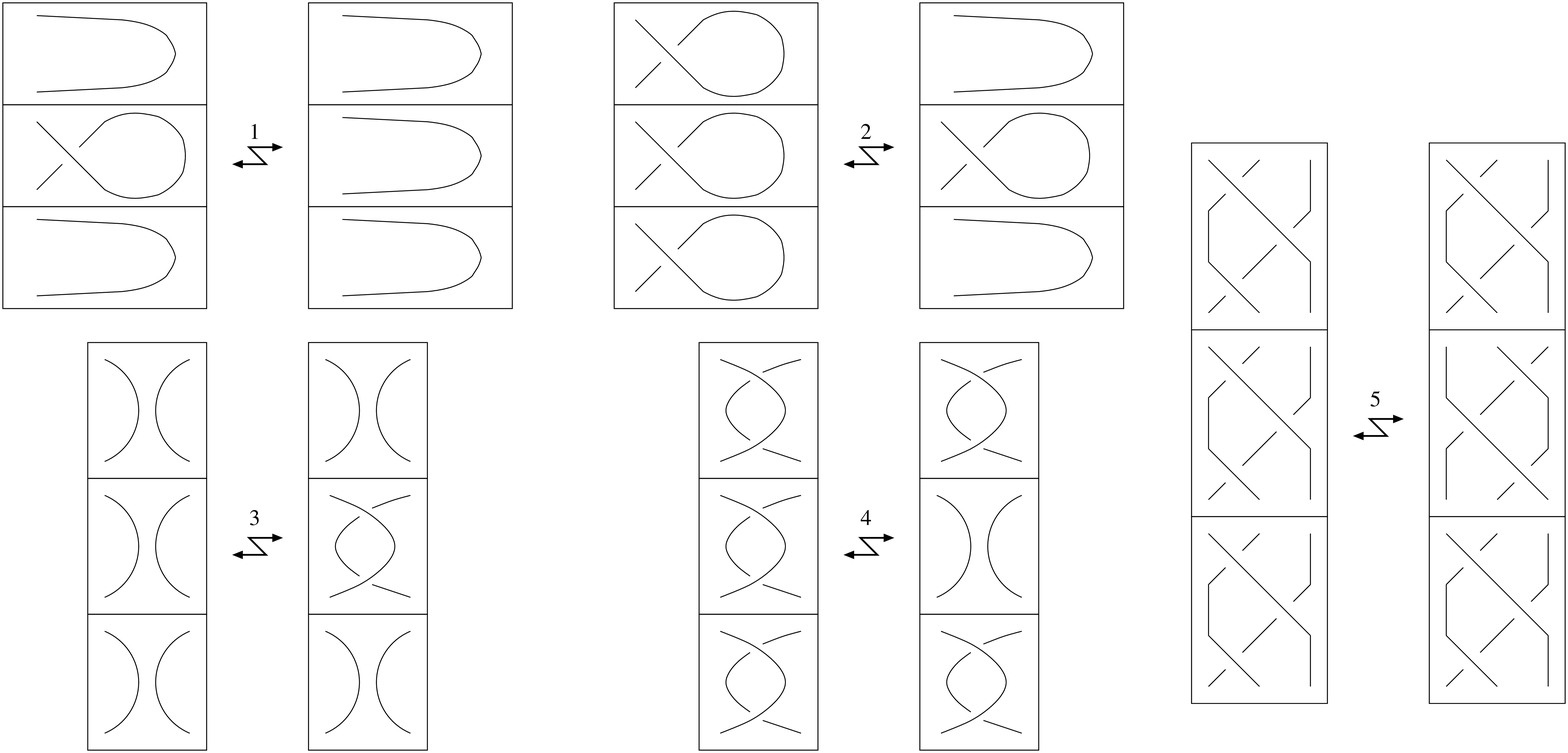}
\caption{Movie moves (Roseman moves)}
\label{1movie5}
\end{center}
\end{figure}

\begin{figure}[ht!]\small
\begin{center}
\psfrag{6}{6}
\psfrag{7}{7}
\includegraphics[width = 11 cm, height = 8 cm]{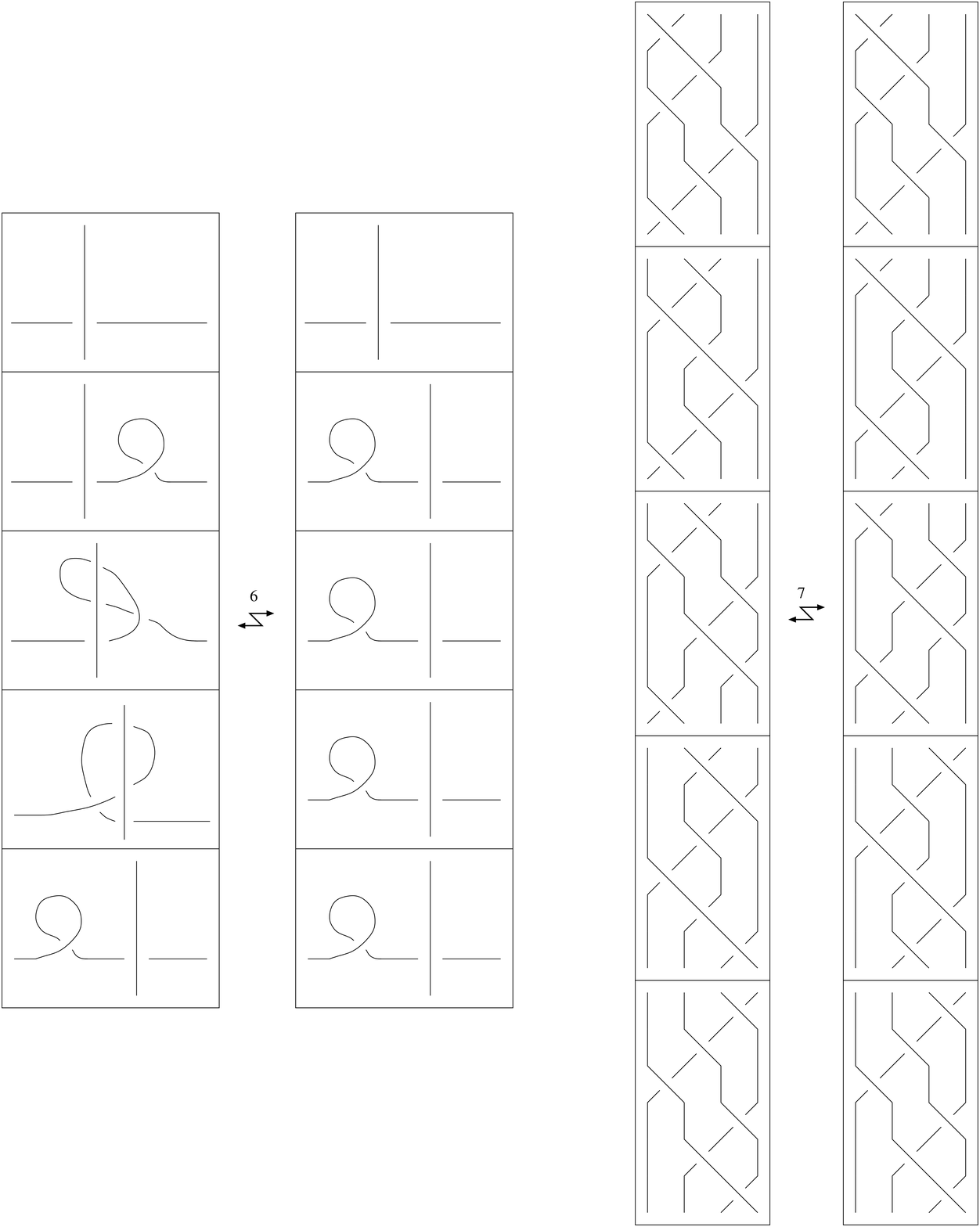}
\caption{Movie moves (Roseman moves)}
\label{6movie7}
\end{center}
\end{figure}

\begin{figure}[ht!]\small
\begin{center}
\psfrag{8}{8}
\psfrag{9}{9}
\psfrag{10}{10}
\psfrag{11}{11}
\psfrag{12}{12}
\psfrag{13}{13}
\psfrag{14}{14}
\psfrag{15}{15}
\includegraphics[width = 12 cm, height = 15 cm]{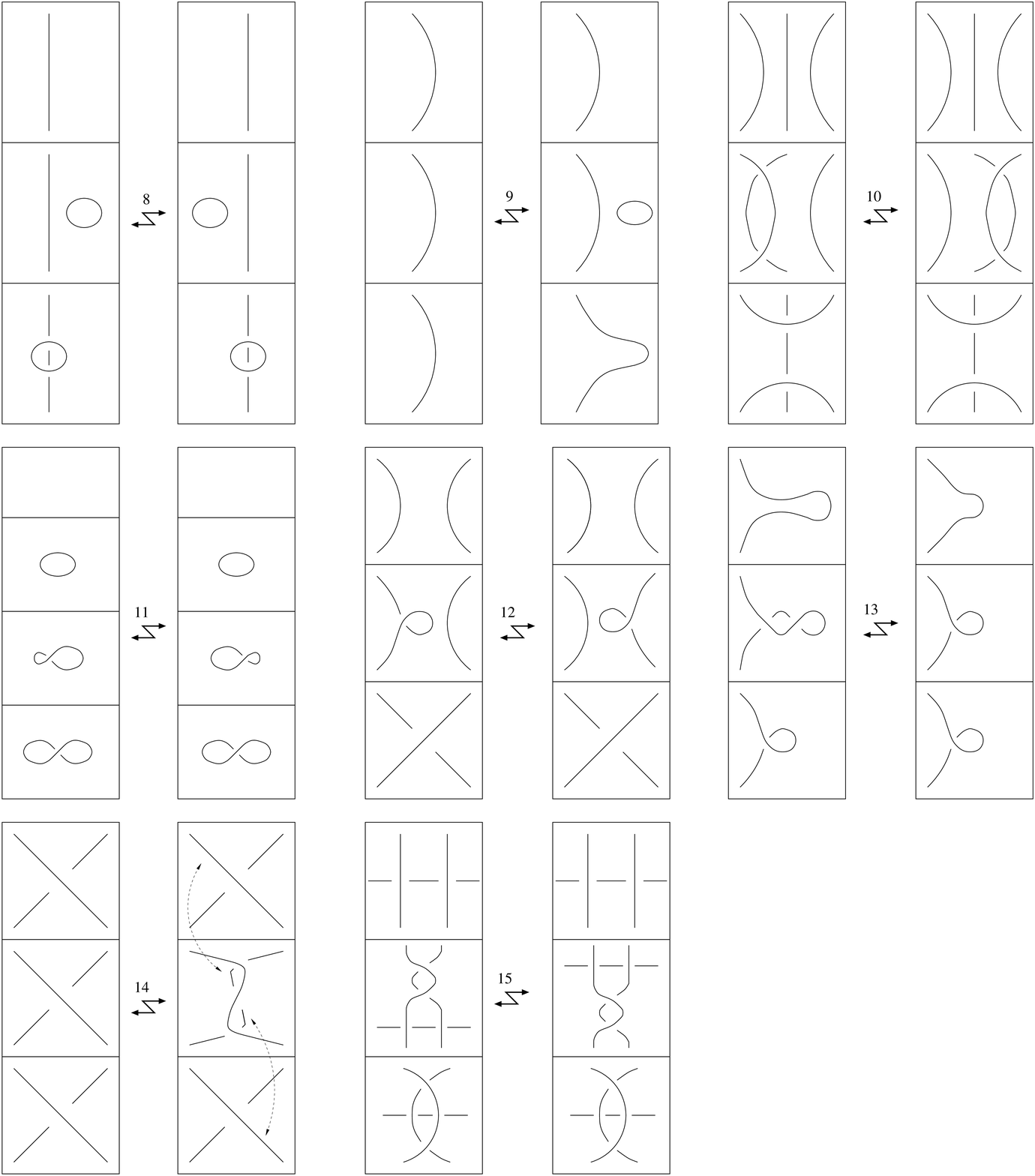}
\caption{Movie moves (additional moves due to Carter/Saito)}
\label{8movie15}
\end{center}
\end{figure}

Carter and Saito \cite{CS} have found the Reidemeister-type moves
for movies. These include movie versions of the Roseman moves mentioned above, 
but also additional moves to handle the additional structure 
of a time function on the diagram. The additional moves do not 
change the local topology of the surface diagram, but the Roseman 
moves do. We give pictures of the {\em movie moves} 
in Figures \ref{1movie5}, \ref{6movie7} and \ref{8movie15}. Each move has two 
sides, which will be called the left and right side according to
their position in these figures. Any two movies of ambient 
isotopic link cobordisms can be related by a sequence of movie 
moves and interchanges of distant critical points.

\begin{remark} There are in fact more movie moves than are displayed 
in the figures. First, each move can be read in two directions, upwards  
or downwards. Second, each move comes in several versions, with varying 
crossing information. Each move has a mirror image, obtained by 
changing all the crossings. Moves 6 and 15 also have additional 
versions depending on the relative positions of the ingoing 
strands. The complicated Move 7 comes in several versions, 
but as we explain in the proof of the main theorem 
(Section \ref{pmkc}), they are all equivalent modulo 
other moves. For details on the movie moves, see for example \cite{CRS}, 
where an even more refined set is described, which is not needed 
for our purposes. Finally, let us remark that it will 
not be necessary to consider any orientations on the strings. 
\end{remark}

\subsection{Local moves}
\label{rmm}
Khovanov associated to each local move on a link diagram 
a homomorphism between the homology groups of its source and 
target diagrams. In this section we describe how this is done. 

\subsubsection{Reidemeister moves}
\label{rm}
Let $D$ and $D'$ be link diagrams which differ by a single Reidemeister 
move of first or second type and assume that $D'$ is the one with more 
crossings. Khovanov proved that the chain complex $C(D')$ splits as 
$C' \oplus C'_{contr}$, where there is an isomorphism 
$\psi: C' \xrightarrow{\cong} C(D)$ and $C'_{contr}$ is chain contractible. 

Thus, the composition $\Psi$ of this isomorphism with the projection 
onto the first summand 
\begin{displaymath}
C(D') = C' \oplus C'_{contr} \xrightarrow{pr_1} C'  \xrightarrow{\psi} C(D)
\end{displaymath}
is a chain equivalence with chain inverse $\Phi$ given by composition of 
the inclusion $i_1$ with $\psi^{-1}$. 

\begin{figure}[ht!]\small
\begin{center}
\psfrag{p}{$p$}
\psfrag{tensorxa}{$\otimes [xa]$}
\psfrag{tensorx}{$\otimes [x]$}
\psfrag{tecken}{$(-1)^i$}
\psfrag{-}{$-$}
\includegraphics[width = 10 cm, height = 1.5 cm]{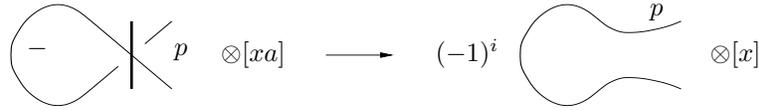}
\caption{The isomorphism $\psi$ for the first Reidemeister move,
  negative twist: the expressions on the left are generators for
  $C'$. The crossing in the twist is called $a$.}
\label{psi1a}
\end{center}
\end{figure}

\begin{figure}[ht!]\small
\begin{center}
\psfrag{tensorx}{$\otimes [x]$}
\psfrag{+}{$+$}
\psfrag{-}{$-$}
\includegraphics[width = 10 cm, height = 3 cm]{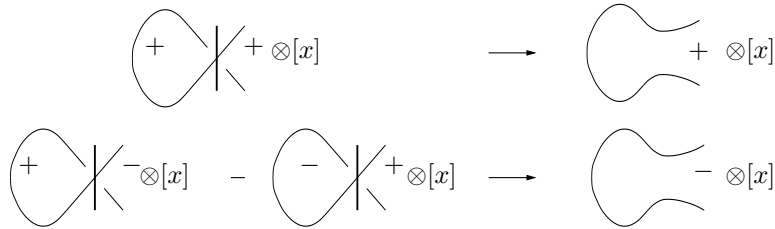}
\caption{The isomorphism $\psi$ for the first Reidemeister move,
  positive twist: the expressions on the left are generators for $C'$. }
\label{psi1b}
\end{center}
\end{figure}

\begin{figure}[ht!]\small
\begin{center}
\psfrag{-}{$-$}
\psfrag{+}{$+$}
\psfrag{p}{$p$}
\psfrag{q}{$q$}
\psfrag{p:q}{$p{:}q$}
\psfrag{q:p}{$q{:}p$}
\psfrag{tecken}{$(-1)^i$}
\psfrag{tensorxa}{$\otimes [xa]$}
\psfrag{tensorxb}{$\otimes [xb]$}
\psfrag{tensorx}{$\otimes [x]$}
\includegraphics[width = 10 cm, height = 2.5 cm]{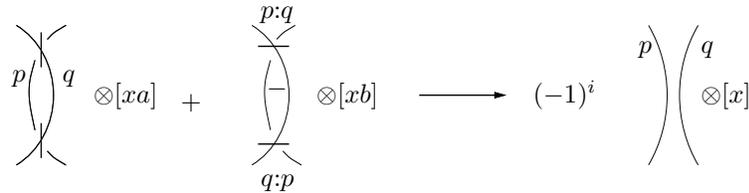}
\caption{The isomorphism $\psi$ for the second Reidemeister move: 
the expressions on the left are generators of $C'$.
The upper crossing is called $a$, the lower $b$.}
\label{psi2}
\end{center}
\end{figure}

\begin{figure}[ht!]\small
\begin{center}
\psfrag{+}{$+$}
\psfrag{-}{$-$}
\psfrag{r}{r}
\psfrag{p}{$p$}
\psfrag{q}{$q$}
\psfrag{a}{$a$}
\psfrag{b}{$b$}
\psfrag{c}{$c$}
\psfrag{p:r}{$p{:}r$}
\psfrag{r:p}{$r{:}p$}
\psfrag{p:q}{$p{:}q$}
\psfrag{q:p}{$q{:}p$}
\psfrag{tensorxa}{$\otimes [xa]$}
\psfrag{tensorxb}{$\otimes [xb]$}
\psfrag{tensorx}{$\phantom{,}\otimes [x]$}
\includegraphics[width = 10 cm, height = 6 cm]{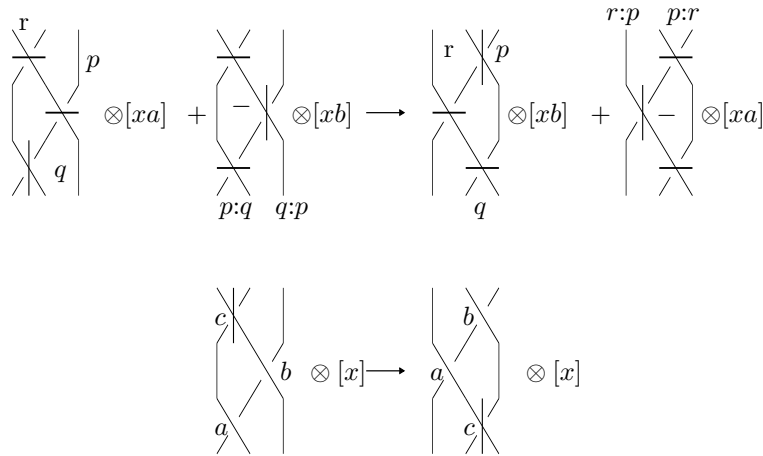}
\caption{The isomorphism $\psi$ for the third Reidemeister
  move: expressions on the left and right are generators for $C$ and
  $C'$ respectively.}
\label{psi3}
\end{center}
\end{figure}

If $D$ differs from $D'$ by a move of the third type, then $C(D)$
and $C(D')$ both split as above, $C(D) = C \oplus C_{contr}$ 
and $C(D') = C' \oplus C'_{contr}$, where there is an isomorphism 
$\psi: C \xrightarrow{\cong}  C'$ and $C_{contr}$ and $C'_{contr}$ are chain 
contractible. Thus these chain complexes are chain equivalent as well
via the map $\Psi$ given by the composition
\begin{displaymath}
C(D') = C' \oplus C'_{contr} \xrightarrow{pr_1} C'  \xrightarrow{\psi}
C \xrightarrow{i_1} C \oplus C_{contr} = C(D).
\end{displaymath}
Generators of the splitting factors, together with the isomorphisms 
$\psi$ are given in Figures \ref{psi1a} through \ref{psi3}.

The figures should be interpreted in the following way.
In the left columns the intersections of the generators 
with the changing disc are displayed. The image of each generator 
under $\psi$ is displayed in the right column. Any state appearing as a 
term in the right column differs from the corresponding state(s) 
in the left column only on markers and circles intersecting the 
changing disc.

An explanation is probably needed for Figure \ref{psi3}. In the first
row, the second term on the left hand side has no letter assigned to its upper
left arc. It is supposed to have the same sign as in the first term,
that is $r$, unless it is connected to one of the other arcs in
the picture. The same applies to the right hand side.

In the second row, no marker is displayed at crossings $a$,$b$ on any
side. We mean that any markers are allowed at those crossings on the
left and that the markers are the same on the right at the 
corresponding crossings. It follows that the resolutions are isotopic, 
and it is understood that the signs are preserved, too.

\begin{figure}[ht!]\small
\begin{center}
\psfrag{p}{$p$}
\psfrag{+}{$+$}
\psfrag{-}{$-$}
\psfrag{p:q}{$p{:}q$}
\psfrag{q:p}{$q{:}p$}
\includegraphics[width = 5 cm, height = 3 cm]{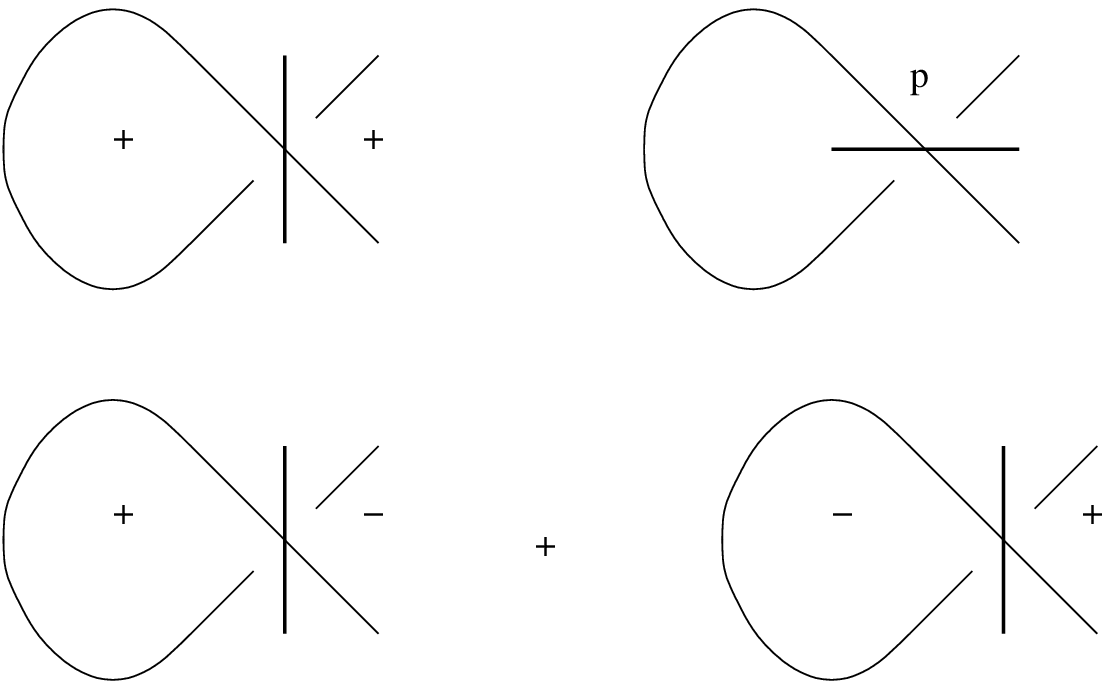}
\caption{States with these configurations in the changing disc generate the 
contractible subcomplex $C_{contr}$ in the case of a negative twist.}
\label{contr1a}
\end{center}
\end{figure}

\begin{figure}[ht!]\small
\begin{center}
\psfrag{p}{$p$}
\psfrag{-}{$-$}
\includegraphics[width = 5 cm, height = 1.5 cm]{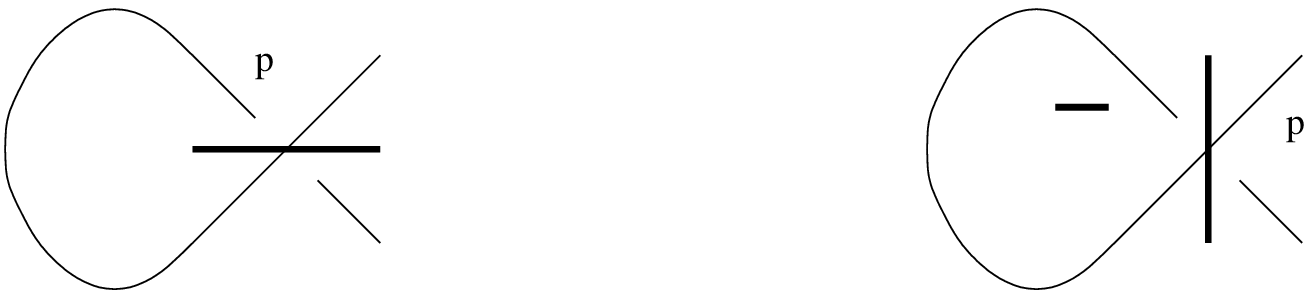}
\caption{States with these configurations in the changing disc generate 
the contractible subcomplex $C_{contr}$ in the case of a positive twist.}
\label{contr1b}
\end{center}
\end{figure}

\begin{figure}[ht!]\small
\begin{center}
\psfrag{p}{$p$}
\psfrag{+}{$+$}
\psfrag{-}{$-$}
\psfrag{p:q}{$p{:}q$}
\psfrag{q:p}{$q{:}p$}
\includegraphics[width = 10 cm, height = 2 cm]{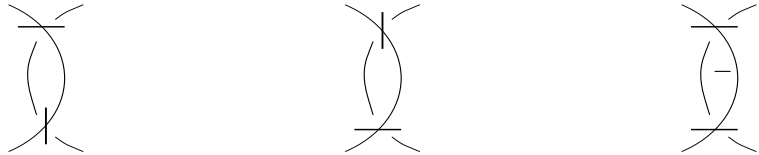}
\caption{Generators for the contractible subcomplex $C_{contr}$, for 
the second Reidemeister move}
\label{contr2}
\end{center}
\end{figure}

\begin{figure}[ht!]\small
\begin{center}
\psfrag{-}{$-$}
\includegraphics[width = 8 cm, height = 5 cm]{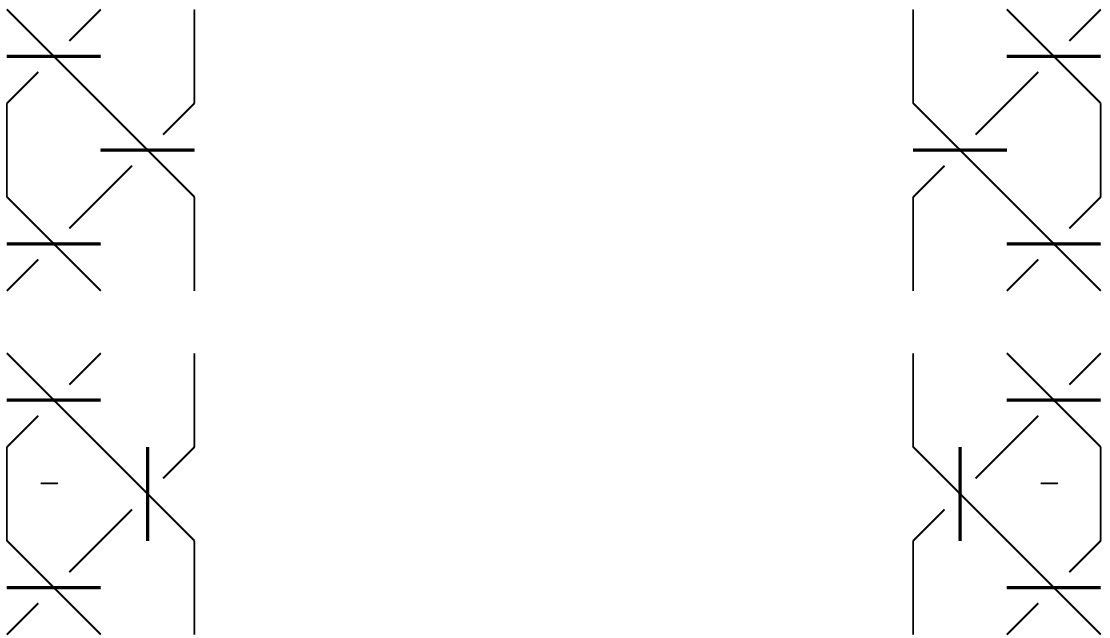}
\caption{Generators for the contractible subcomplexes $C_{contr}$ 
and $C'_{contr}$, for the third Reidemeister move:   note that they 
coincide after a rotation by $\pi$.}
\label{contr3}
\end{center}
\end{figure}

Generators of the contractible summands are given in Figures \ref{contr1a} 
through \ref{contr3}. 

\subsubsection{Morse moves}
\begin{figure}[ht!]\small
\begin{center}
\psfrag{+}{$+$}
\psfrag{-}{$-$}
\psfrag{p}{$p$}
\psfrag{q}{$q$}
\psfrag{p:q}{$p{:}q$}
\psfrag{q:p}{$q{:}p$}
\psfrag{minimum}{minimum}
\psfrag{maximum}{maximum}
\psfrag{saddle point}{saddle point}
\includegraphics[width = 7 cm, height = 6 cm]{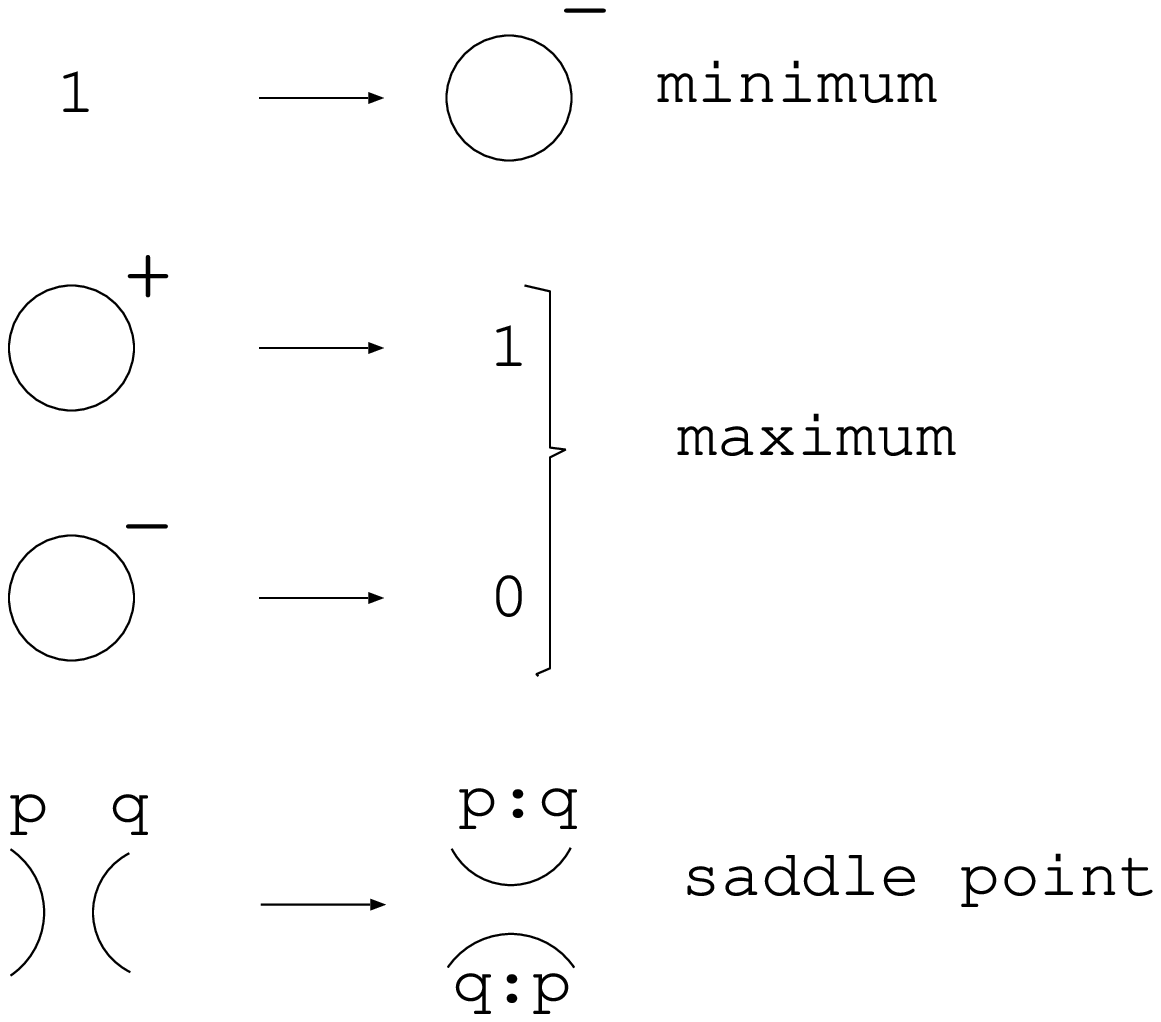}
\caption{The effect of Morse modifications on states}
\label{chain4}
\end{center}
\end{figure}

Khovanov also introduced a chain map of bidegree $(0, 1)$ between 
the complexes of two diagrams which differ by a single Morse 
modification corresponding to a maximum or minimum, 
and a chain map of bidegree $(0,-1)$ if they differ by a 
saddle-point modification. These maps are described by Figure \ref{chain4}.

\subsubsection{Chain maps in the basis of states}
\begin{figure}[ht!]\small
\begin{center}
\psfrag{tecken}{$(-1)^i$}
\psfrag{p}{$p$}
\psfrag{tecken+1}{$(-1)^{i+1}$}
\psfrag{tensorxa}{$\otimes [xa]$}
\psfrag{tensorx}{$\otimes [x]$}
\psfrag{+}{$+$}
\psfrag{-}{$- $}
\psfrag{Psi}{$\Psi$}
\psfrag{Phi}{$\Phi$}
\includegraphics[width = 10 cm, height = 10 cm]{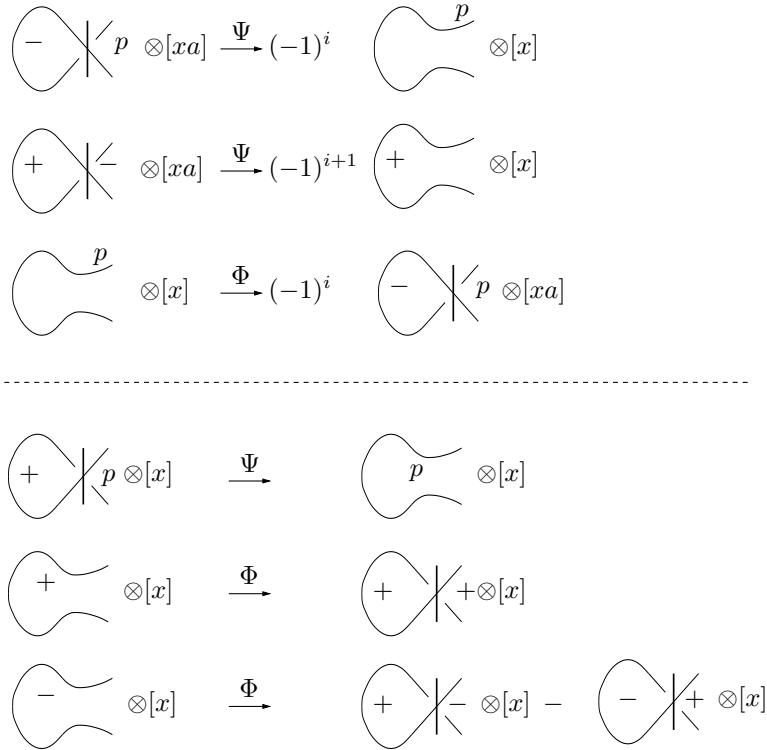}
\caption{The effect of the first Reidemeister moves on states:
the negative twist is above the dashed line, the positive twist below.
States with any other local configuration map to zero.}
\label{chain1}
\end{center}
\end{figure}

\begin{figure}[ht!]\small
\begin{center}
\psfrag{p}{$p$}
\psfrag{q}{$q$}
\psfrag{p:q}{$p{:}q$}
\psfrag{q:p}{$q{:}p$}
\psfrag{tkparent}{$(-1)^i$ (}
\psfrag{tecken+1}{$(-1)^{i+1}$}
\psfrag{tecken}{$(-1)^i$}
\psfrag{tensorxa}{$\otimes [xa]$}
\psfrag{tensorx}{$\otimes [x]$}
\psfrag{xbparent}{$\otimes [xb]$ )}
\psfrag{tensorxb}{$\otimes [xb]$}
\psfrag{+}{$+$}
\psfrag{-}{$-$}
\psfrag{Psi}{$\Psi$}
\psfrag{Phi}{$\Phi$}
\includegraphics[width = 12 cm, height = 8 cm]{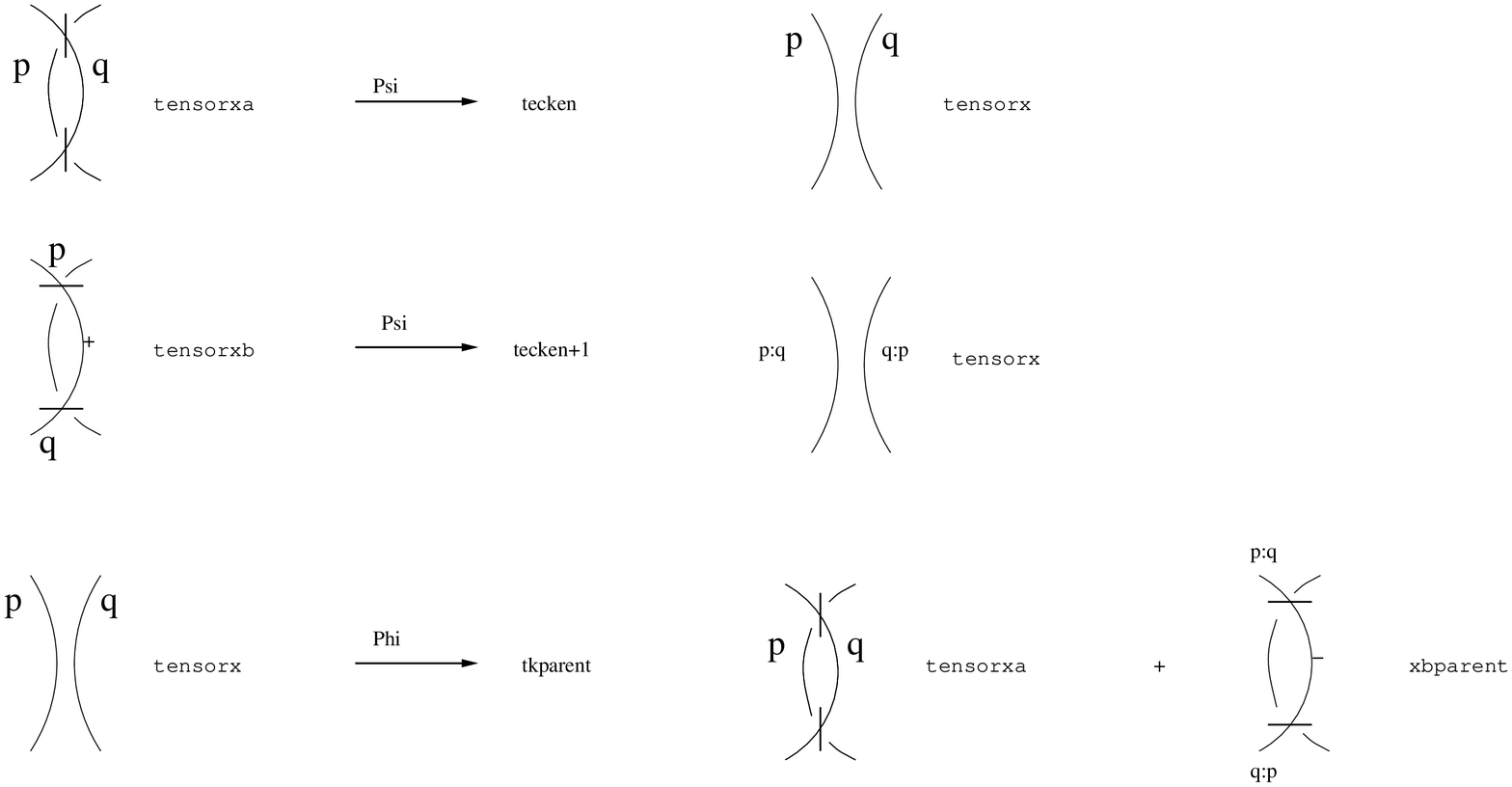}
\caption{The effect of the second Reidemeister move on states:
states with any other local configuration map to zero.}
\label{chain2}
\end{center}
\end{figure}

\begin{figure}[ht!]\small
\begin{center}
\psfrag{P}{$\Psi$}
\psfrag{p}{$p$}
\psfrag{q}{$q$}
\psfrag{r}{$r$}
\psfrag{+}{$+$}
\psfrag{-}{$-$}
\psfrag{p:q}{$p{:}q$}
\psfrag{q:p}{$q{:}p$}
\psfrag{p:r}{$p{:}r$}
\psfrag{r:p}{$r{:}p$}
\psfrag{r:(q:p)}{$r{:}(q{:}p)$}
\psfrag{(q:p):r}{$(q{:}p){:}r$}
\psfrag{tensorxa}{$\otimes [xa]$}
\psfrag{tensorxb}{$\otimes [xb]$}
\psfrag{tensorxab}{$\otimes [xab]$}
\psfrag{tensorxc}{$\otimes [xc]$}
\psfrag{tensorxbc}{$\otimes [xbc]$}
\psfrag{tensorx}{$\otimes [x]$}
\includegraphics[width = 13 cm, height = 9 cm]{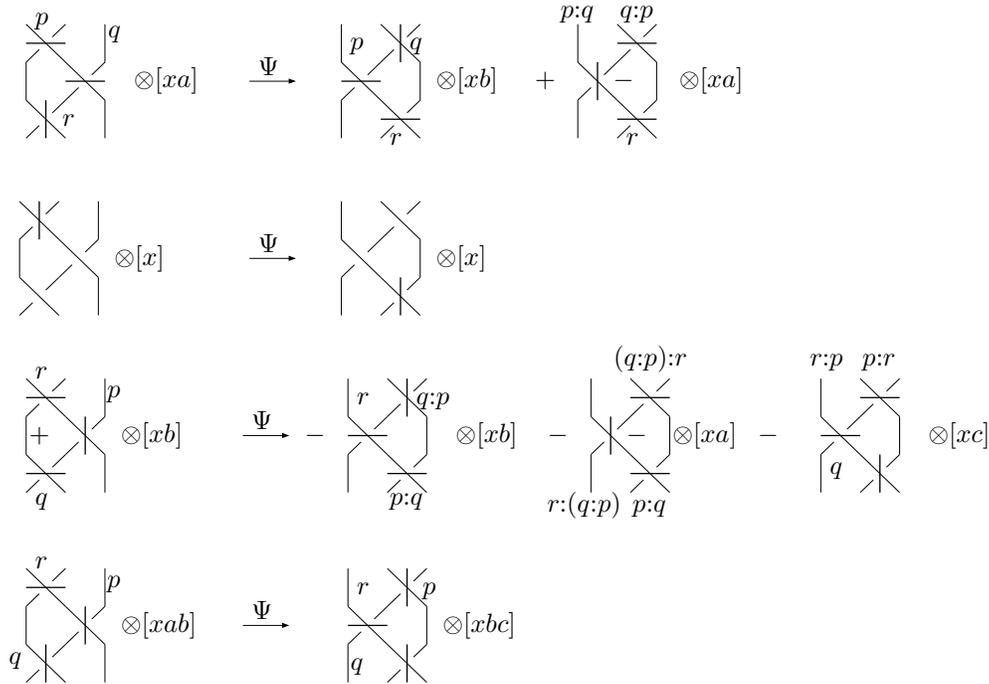}
\caption{The effect of the third Reidemeister move on states: 
states with any other local configuration map to zero. The 
crossings are, from bottom to top, a,b,c on the left and 
c,a,b on the right.}
\label{chain3}
\end{center}
\end{figure}

The states form a canonical basis (up to reversing signs) 
for Khovanov's chain complex. In the following proposition 
we give an explicit description of the value of the chain 
maps on the states. Recall the chain equivalences $\Psi$  
and their homotopy inverses $\Phi$ from Section \ref{rm}. 
\begin{proposition} 
\label{prop}
The values of the chain maps $\Psi$ and $\Phi$ on states are as shown
in Figure \ref{chain1} and \ref{chain2} for the first two 
Reidemeister moves and in Figure \ref{chain3} for the third 
Reidemeister move. The Morse move maps are as shown in Figure \ref{chain4}. 
\end{proposition}

\begin{proof} For Morse moves there is nothing to prove.
We rewrite the chain equivalences $\Psi$ and $\Phi$
in the basis of states.  We give a sample calculation here 
and refer to \cite{Ja1} for more details. 

Let us explain the second row in Figure \ref{chain2}. 
We call the upper crossing $a$ and the lower crossing $b$.

Let $S_{+-} \otimes [xb]$ be the state on the left of this row.
We use the index ``{$+-$}'' to stand for ``positive marker at $a$ 
and negative at $b$''. 
Then there is a unique state $S_{++} \otimes [x]$, such that 
$S_{+-}$ and $S_{++}$ coincide outside the changing
disc. Note that $S_{++} \otimes [x]$ lies in $C_{contr}$ 
(see Figure \ref{contr2}).
Now apply the differential to this state.
\begin{displaymath}
d(S_{++} \otimes [x]) = S_{+-} \otimes [xb] + S_{-+} \otimes [xa] + 
(S_{+-,-} \otimes [xb]) + \sum_t S_{++}^t \otimes [xt]
\end{displaymath}
Here, $t$ ranges over those crossings outside the changing disc 
where $S_{++}$ has a positive marker. $S_{+-,-}$ has positive marker 
at $a$, negative marker at $b$ and a negative sign on the small circle
enclosed in the changing disc. It appears only if the sign of the
bottom arc in $S_{+-}$ is negative.

The left hand side of this equation is in $C_{contr}$ since $S_{++}$ is.
So is the third term and the big sum on the right hand side.
Modulo the contractible subcomplex $C_{contr}$, we are left with
\begin{displaymath}
0 = S_{+-} \otimes [xb] + S_{-+} \otimes [xa].
\end{displaymath}
By the definition of the differential the signs of the circles of
$S_{-+}$ depend on those of $S_{+-}$ as determined by the ``Frobenius 
calculus of signed circles'' in Figure \ref{frobenius}. 
Applying the isomorphism $\psi$ 
to this expression it follows that the value on $S$ is as 
displayed in Figure \ref{chain2}. 
\end{proof}

\begin{figure}[ht!]\small
\begin{center}
\psfrag{o32}{$\overline{\Omega}_3$}
\psfrag{o3}{$\Omega_3$}
\psfrag{o2}{$\Omega_2$}
\psfrag{isotopy}{isotopy}
\includegraphics[width = 7 cm, height = 7 cm]{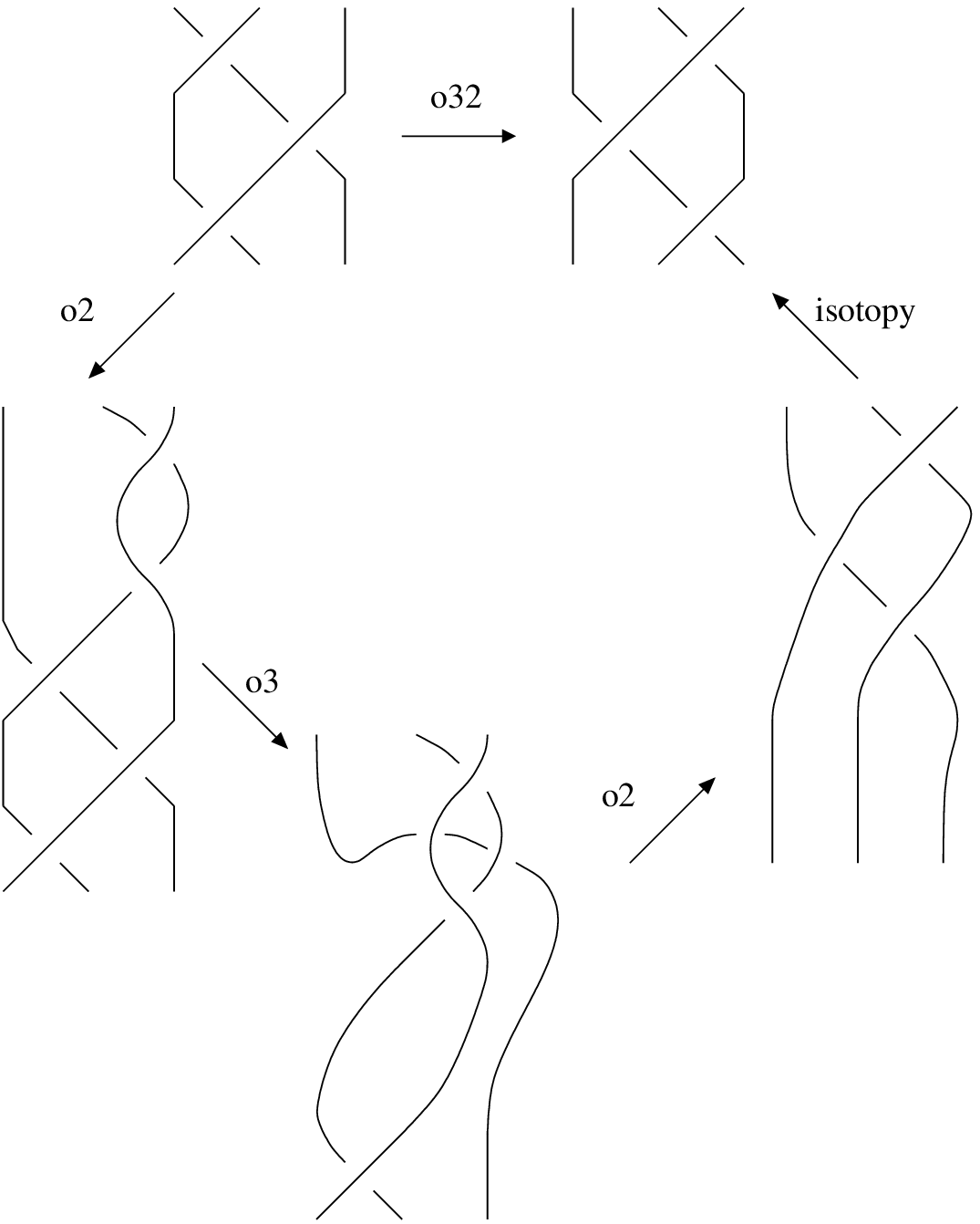}
\caption{Another third Reidemeister move}
\label{turaev}
\end{center}
\end{figure}

\subsection{Yet another move}
\label{yetanothermove}
Below we use the notation $\Omega_k$ for the $k$-th Reidemeister move.
We will need one additional move whose induced map is not yet 
explicitly described; the mirror image $\overline{\Omega}_3$ of $\Omega_3$. 
The chain equivalence corresponding to this move was not 
mentioned by Khovanov in \cite{K} since it was not needed for his purposes
there. It is well-known that it can be expressed in the 
other moves. We will fix one such expression to ensure  
that we have a unique induced homomorphism. 
To this end, consider Figure \ref{turaev}. 
There $\overline{\Omega}_3$ is expressed as the composition of an 
$\Omega_2$-move, the original $\Omega_3$-move
and a second $\Omega_2$-move. 
We have a choice as to where to make the first $\Omega_2$.
We choose, as shown in the picture, the wedge formed by 
the upper and middle strands. The composition induces 
a well-defined isomorphism on homology, in view of the 
previous sections.


\subsection{Definition of the invariant}
\label{def}
Let $(\Sigma, L_0, L_1)$ be a link cobordism and let a movie
presentation of $\Sigma$ be given. 
Close to a critical level, the movie realizes 
a single Reidemeister move or oriented Morse modification.
It follows that each critical level induces a chain map 
between the complexes of the diagrams just before and 
after this level. Between two consecutive critical levels 
the diagram undergoes a planar isotopy. There is an obvious 
canonical isomorphism associated to this planar isotopy, 
given by tracing each state through it.

Composing all these maps, we obtain a chain homomorphism 
$\phi_{\Sigma}$ from $C^{i,j}(D_0)$ to $C^{i,j + \chi (\Sigma)}(D_1)$.
 
The reason for the grading of this map is plain. 
The Reidemeister moves are grading preserving, and the Morse 
modifications change the grading by $+1$ for extrema and by $-1$ for
saddle points. That is, for each attachment of a $0$-- or $2$--handle 
the grading increases by $1$ and for each attachment of a $1$--handle 
the grading decreases by $1$. The sum of these changes is the 
Euler characteristic of $\Sigma$.

The chain map $\phi_{\Sigma}$ induces a map $\phi_{\Sigma *}$ 
on homology and the purported invariant of link cobordisms is defined by 
\begin{displaymath}
\Sigma \mapsto \pm \phi_{\Sigma *}.
\end{displaymath}
The invariance is shown in Section \ref{pmkc}. 

\begin{remark} A closed knotted surface $\Sigma$ in 4--space is a cobordism 
between empty links. The empty link has a single homology group 
$\mathbb{Z}$ in degree $(0,0)$, so in this case the invariant is a single
natural number. If the Euler characteristic of $\Sigma$ is non-zero 
this number must be zero, because of the above grading properties. 
Finally note that for an unknotted torus the number is $2$.
\end{remark}


\section{Khovanov's conjecture}
\label{ce}
\subsection{Original statement and examples}
\label{cekc}
In \cite{K}, the following conjecture was made.

\begin{conjecture}(Khovanov, \cite{K}, page 40)\qua If two movies of a 
link cobordism $\Sigma$ have the property that the corresponding 
source and target link diagrams of the two movies are isomorphic, 
the induced map on homology is the same up to an overall sign. 
\end{conjecture}

Any ambient isotopy of a link in 3--space gives rise to a link
cobordism, which is ambient isotopic to a cylinder on the link,
i.e.\ to the trivial cobordism. Furthermore, the effect of an 
ambient isotopy on a link diagram is a sequence of planar isotopies 
and Reidemeister moves. Hence, if taken literally, this conjecture implies 
that any sequence of planar isotopies and Reidemeister moves, 
which starts and finishes with the same link diagram either 
induces the identity map in all of the homology groups 
associated to this diagram, or induces multiplication by 
$-1$ in all of them. This is not true, even for very simple 
link diagrams, as the following example shows.

\begin{remark} It is obvious that the homology groups of the 
unnested unlink diagram ${0}_1^2$ of two components without crossings 
coincide with the chain groups. 
\begin{displaymath}
\begin{split}
C^{0,-2} = \, &\mathcal{H}^{0,-2} \cong \mathbb{Z} \\
C^{0,0} = \,&\mathcal{H}^{0,0} \,\, \cong \mathbb{Z} \oplus \mathbb{Z} \\
C^{0,2} = \,&\mathcal{H}^{0,2} \,\, \cong \mathbb{Z} \\
\end{split}
\end{displaymath}
A planar isotopy which interchanges the two components also
interchanges the two generators of $C^{0,0}$ (the states which 
have different signs on the two circles) and is neither 
plus or minus the identity. (This fact was first pointed 
out by Olof-Petter \"Ostlund \cite{O}.) Clearly, similar things 
can be done with more and non-trivial components. The 
question whether non-trivial maps can be induced 
on the homology of a {\em knot} diagram is now natural. 
In Theorem \ref{818} below we find an affirmative answer to that question.
\end{remark}

\begin{definition}
The subgroup of the automorphism group of $\mathcal{H}(D)$ consisting of
isomorphisms induced from sequences of Reidemeister moves and planar
isotopies from $D$ to itself we call {\em the monodromy group} of $D$.
The monodromy group is called {\em trivial} if it consists of 
only the identity map and its negative.
\end{definition}

\begin{remark} The monodromy group of $0_1^2$ above is non-trivial.
\end{remark}

\begin{figure}[ht!]\small
\begin{center}
\includegraphics[width = 3 cm, height = 3 cm]{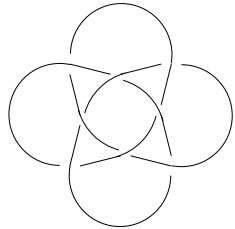}
\caption{The knot diagram $8_{18}$}
\label{knot818}
\end{center}
\end{figure}

\begin{theorem}
\label{818}
The monodromy group of $\mathcal{H}^{i,j}(D)$ is non-trivial when $D$ is the
diagram of the knot $8_{18}$ pictured in Figure \ref{knot818}. 
\end{theorem}

\begin{proof}
We will compute $C^{i,j}(D)$ for $j=-7$. Then the following formulas
obviously hold. 
\begin{displaymath}
\begin{split}
&w(D) = 0 \\
&\sigma(s) + 2\tau(S) = 14 \\
&-8 \leq \sigma(s) \leq 8 \\
\end{split}
\end{displaymath}
If all markers in $S$ are positive then $\sigma(S) = 8$, i.e.\ $i = -4$. 
The resolution of $S$ consists of five circles, four of which are
pairwise unnested and enclosed by the fifth. By the second of the above
equations $\tau(S)=3$, so one of the circles has to be negative and
the rest positive. There are five such states, $v_1,...,v_5$, and they 
are shown in Figure \ref{vs}. 

\begin{figure}[ht!]\small
\begin{center}
\includegraphics[width = 13 cm, height = 6 cm]{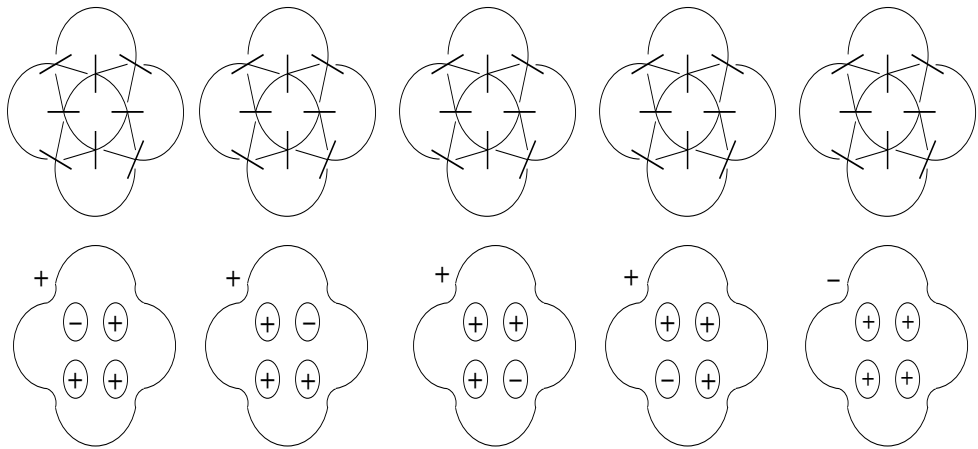}
\caption{Generators $v_1,...,v_5$ of $C^{-4,-7}(8_{18})$}
\label{vs}
\end{center}
\end{figure}

If $i=-3$, then $\sigma = 6$ so $\tau = 4$. 
The resolution of any such state (i.e.\ with exactly one negative
marker) consists of three unnested circles enclosed by a fourth. 
All circles are equipped with positive signs. There are eight states of 
this form, $u_1,...,u_4$ and $w_1,...,w_4$. They are shown in
Figure \ref{uws}. 

\begin{figure}[ht!]\small
\begin{center}
\includegraphics[width = 11 cm, height = 11 cm]{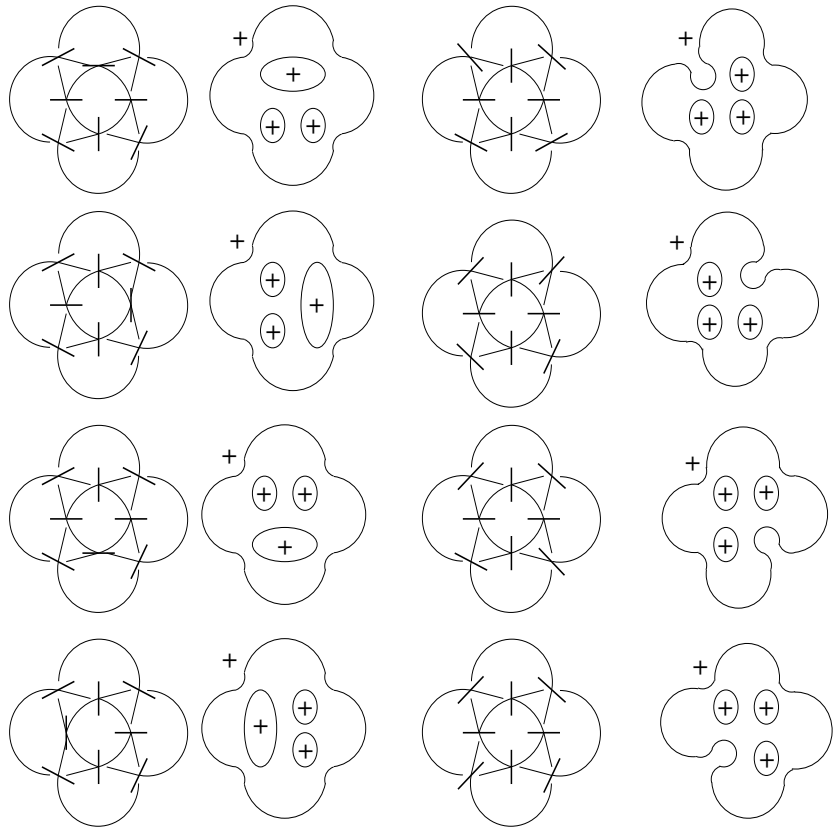}
\caption{Generators $u_1,...,u_4$ (left column) and $w_1,...,w_4$
(right column) of $C^{-3,-7}(8_{18})$}
\label{uws}
\end{center}
\end{figure}

The differential of any state $S$ in $C^{-3,-7}(D)$ is zero since any change of
a positive marker in $S$ causes two positive circles to merge. (In
fact, all groups $C^{i,-7}(D)$ are zero for $i$ different from $-3$ or $-4$.) 
We have the following fragment of the chain complex:
\begin{displaymath}
\xymatrix{
0 \ar[r] & C^{-4,-7}(D) \ar[d]^\cong \ar[r]^d & C^{-3,-7}(D)
\ar[d]^\cong \ar[r] & 0 \\
0 \ar[r] & \mathbb{Z}^5 \ar[r]^A & \mathbb{Z}^8 \ar[r] & 0}
\end{displaymath}
where the vertical isomorphisms are given by the choices of ordered bases
$v_i \mapsto e_i$ and $w_i \mapsto e_i$, $u_i \mapsto e_{i+4}$, where
$\{e_i\}_{i=1}^n$ is the canonical basis in $\mathbb{Z}^n$. 
The definition of $d$ gives:
\begin{displaymath}
\begin{split}
&v_i \mapsto w_i + u_i + u_{i-1},\quad 1 \leq i \leq 4\\
&v_5 \mapsto \sum_i w_i\\
\end{split}
\end{displaymath}
with the natural cyclic ordering of indices. We do not need to
consider the $E(L)$-tensor factor, since for any state 
involved its singleton or empty set of negative markers is uniquely 
ordered. Thus, the matrix $A$ is given by 
\[
A = 
\left(
\begin{array}{ccccc}
1 & 0 & 0 & 0 & 1 \\
0 & 1 & 0 & 0 & 1 \\
0 & 0 & 1 & 0 & 1 \\
0 & 0 & 0 & 1 & 1 \\
1 & 1 & 0 & 0 & 0 \\
0 & 1 & 1 & 0 & 0 \\
0 & 0 & 1 & 1 & 0 \\
1 & 0 & 0 & 1 & 0 \\
\end{array}
\right)
\sim
\left(
\begin{array}{ccccc}
1&0&0&0&0\\
0&1&0&0&0\\
0&0&1&0&0\\
0&0&0&1&0\\
0&0&0&0&2\\
0&0&0&0&0\\
0&0&0&0&0\\
0&0&0&0&0\\
\end{array}
\right)
= A'
\]
and $A'$ is the Smith normal form of $A$. 
It follows that $d$ is injective and that $A$ and $A'$ are
presentation matrices of the isomorphism class of $\mathcal{H}^{-3,-7}(D)$. 
From $A'$ we see that 
$\mathcal{H}^{-3,-7}(D) \cong \mathbb{Z} \oplus \mathbb{Z} \oplus \mathbb{Z}
\oplus \mathbb{Z}_2$. 

Let $\Sigma$ be the planar isotopy which rigidly rotates $D$
clockwise by an angle $\pi / 2$ about the center. 
Denote by $\phi_{\Sigma}$ and $\phi_{\Sigma *}$ the maps induced 
on the chain complex and on homology, respectively.
 
Then it is clear that $\phi_\Sigma(w_i) =
w_{i+1}$ and that $\phi_\Sigma(u_i) =
u_{i+1}$ so that $\phi_{\Sigma *}([w_i])=[w_{i+1}]$ and 
$\phi_{\Sigma *}([u_i])=[u_{i+1}]$. 
Now, suppose the monodromy group were trivial. Then either
$\phi_{\Sigma *}[w_i]=[w_i]$  and $\phi_{\Sigma *}[u_i]=[u_i]$ 
or $\phi_{\Sigma *}[w_i]=-[w_i]$ and $\phi_{\Sigma *}[u_i]=-[u_i]$ . 
This gives two cases:

Case($+$): Then for all $i$, $[w_i] = [w_{i+1}]$ and $[u_i] =
[u_{i+1}]$. This immediately reduces the number of generators 
of $\mathcal{H}^{-3,-7}(D)$ to two and the relations to
\begin{displaymath}
\begin{split}
&w_1 + 2u_1 = 0\\
&4 w_1 = 0 \\
\end{split}
\end{displaymath}
which is a presentation of $\mathbb{Z}_8$. This contradicts the 
computation of $\mathcal{H}^{-3,-7}(D)$ above.

Case($-$): In this case $[w_i]=-[w_{i+1}]$ and $[u_i] = -[u_{i+1}]$. 
This reduces the relations of $\mathcal{H}^{-3,-7}(D)$ to 
\begin{displaymath}
\begin{split}
&w_i = 0,  \quad 1 \leq i \leq 4\\
&u_i = -u_{i+1}, \quad  1 \leq i \leq 4\\
\end{split}
\end{displaymath}
which is a presentation of $\mathbb{Z}$, and we have 
another contradiction. This completes the proof.
\end{proof}

\begin{remark} Let $\Sigma$ be the cobordism of the planar isotopy 
in the proof above. With coefficients in $\mathbb{Q}$, 
$\mathcal{H}^{-3,-7}(D) \cong \mathbb{Q}^3$. Three independent cycles 
in the ($\phi_{\Sigma}$--invariant) orthogonal complement of 
$dC^{-4,-7}(D)$ are 
\begin{displaymath}
\begin{split}
v_1 &= w_1 - w_3 - u_1 + u_2\\
v_2 &= w_1 + w_2 - w_3 - w_4 - u_1 + u_3\\
v_3 &= w_2 - w_4 - u_1 + u_4\\ 
\end{split}
\end{displaymath}
and it is easy to see that in the basis given by the homology classes 
of these cycles the map $\phi_{\Sigma *}$ on $\mathcal{H}^{-3,-7}(D)$ 
is given by the following matrix.
\[
\left(
\begin{array}{ccc}
-1&-1&-1\\
1&0&0\\
0&1&0\\
\end{array}
\right)
\]
This map is of order four and its spectrum is $\{-1,i,-i\}$.
\end{remark}


\subsection{Recovering the conjecture}
\label{rtc}
Because of the examples in the previous subsection it is 
necessary to reformulate the conjecture to have any hope of 
proving it. It turns out that it is enough to require that the 
isotopy is relative to the boundary.  

\begin{theorem}[Khovanov's conjecture]
\label{mkc}
For oriented links $L_0$ and $L_1$, presented by diagrams $D_0$ and
$D_1$, an oriented link cobordism $\Sigma\subset\mathbb{R}^3\times[0,1]$ 
from $L_0$ to $L_1$, defines homomorphisms $\mathcal{H}^{i,j}(D_0)\to
\mathcal{H}^{i,j+\chi(\Sigma)} (D_1)$ invariant up to an overall 
multiplication by $-1$ under ambient isotopy of $\Sigma$ leaving $\partial
\Sigma$ setwise fixed. Moreover, this invariant is non-trivial.
\end{theorem}
The proof of Theorem \ref{mkc} occupies the next section.
It uses the following result by Carter and Saito, previously 
alluded to in Section \ref{mm}. It is given a slight reformulation 
to fit our purposes.

\begin{theorem*}[Carter--Saito, \cite{CS}] 
\label{cs}Two movies represent 
equivalent link cobordisms in the sense of Theorem \ref{mkc}
if and only if they can be related via a finite sequence of movie moves
and interchanges of distant critical points of the time function $t$.
\end{theorem*}


\section{Proof of Khovanov's conjecture}
\label{pmkc}

We split the proof into three lemmas.

\begin{figure}[ht!]\small
\begin{center}
\includegraphics[width = 5 cm, height = 6.5 cm]{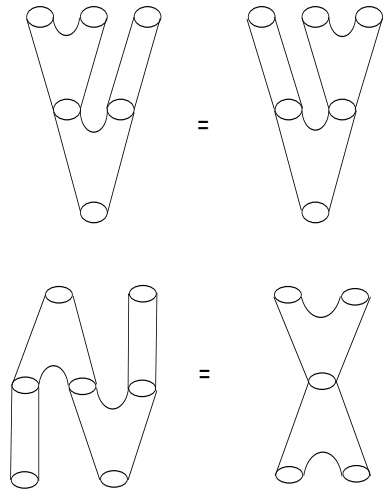}
\caption{Associativity/coassociativity and an additional relation
prove the commutativity of distant critical points.}
\label{frobeniuscommute}
\end{center}
\end{figure}

\begin{lemma}[Distant critical points] The interchange of two distant 
critical points of the surface diagram does not change the induced map
on homology.
\end{lemma}  
\noindent
\begin{proof} Let $D$ be a link diagram and let $D'$ be obtained 
from $D$ by two local moves in disjoint changing discs. 
There are two different orders in which these moves can be performed.
They correspond to two different movies from $D$ to $D'$, inducing
maps $\phi_l, \phi_r: C(D) \rightarrow C(D')$. Let $S$ be a state
of $D$. 

If we temporarily forget the signs of the resolutions the moves are 
completely local and only affect $S$ inside the changing discs. 

From the Figures \ref{chain4} through \ref{chain3} describing 
the induced maps it is clear
that the resolution res($T$) of each term $T$ of $\phi_l(S)$ or 
$\phi_r(S)$ can be obtained from res($S$) by a sequence of
applications of the signed circle calculus (recall Figure
\ref{frobenius}) inside the changing discs, together with the 
addition or removal of some circle enclosed in the
changing disc. The latter changes are local and commute with the others. 
Hence without loss of generality we may forget about them. It is 
therefore sufficient to prove that two distant saddle point moves performed 
on the resolution commute.

This follows from the fact that the product $m$ and coproduct $\Delta$
of the Frobenius algebra $A$ are associative and coassociative
respectively, and satisfy the relation
\begin{displaymath}
\Delta \circ m  = (m \otimes id)(id \otimes \Delta)
\end{displaymath}
In Figure \ref{frobeniuscommute} these relations are interpreted 
topologically. We immediately see that they express the commutativity 
of saddle point moves. This completes the proof.
\end{proof}

\begin{lemma}[Carter--Saito moves]For any Carter--Saito move $M$, 
the maps induced on homology 
by the two sides of $M$ coincide up to an overall sign. 
The sign difference is given in Table \ref{table}.\end{lemma}

\begin{table}[ht!]\small\centerline{
\begin{tabular}{|l|c|c|}
\hline
Movie move no. & Downward time & Upward time \\ 
\hline
\hline
1,2,3,4,5      & same sign        & same sign  \\
\hline
6, pos twist   & different sign   & different sign     \\         
\hline
6, neg twist    & same sign        & same sign  \\ 
\hline
7              & same sign        & same sign  \\
\hline
8              & same sign        & different sign \\ 
\hline
9              & same sign        & same sign  \\
\hline
10             & different sign   & same sign  \\
\hline
11, as displayed   & different sign   & same sign  \\
\hline
11, mirror image  & same sign        & same sign   \\
\hline
12, as displayed   & different sign   & same sign  \\
\hline
12, mirror image  & same sign        & different sign  \\
\hline
13, as displayed   & same sign        & same sign   \\
\hline
13, mirror image   & different sign   & different sign  \\     
\hline
14, as displayed & different sign   & different sign   \\    
\hline
14, mirror image & same sign        & same sign   \\
\hline
15             & different sign   & different sign     \\  
\hline
\end{tabular}}\nocolon\caption{}\label{table}\end{table}

\begin{remark} The words ``Downward'' and ``Upward'' refer to the direction
of the movies in Figures \ref{1movie5}, \ref{6movie7} and \ref{8movie15}. 
``As displayed'' means as displayed in these figures. Observe 
that if the movies contain only Reidemeister moves, the induced maps are 
inverses, and it is enough to consider one time direction. 
\end{remark}

\proof We compute the maps on the chain level and analyze the
effect on homology. 

\subsubsection*{Movie moves 1-5}
The left hand side of each of these moves trivially 
induces the identity. The same holds for the right hand side, since it 
is the composition of a Reidemeister move and its inverse. 

\subsubsection*{Movie move 6} 

\begin{figure}[ht!]\small
\begin{center}
\psfrag{p}{$p$}
\psfrag{q}{$q$}
\psfrag{+}{$+$}
\psfrag{-}{$-$}
\psfrag{p:q}{$p{:}q$}
\psfrag{q:p}{$q{:}p$}
\psfrag{tecken}{$(-1)^i$}
\includegraphics[width = 13 cm, height = 15 cm]{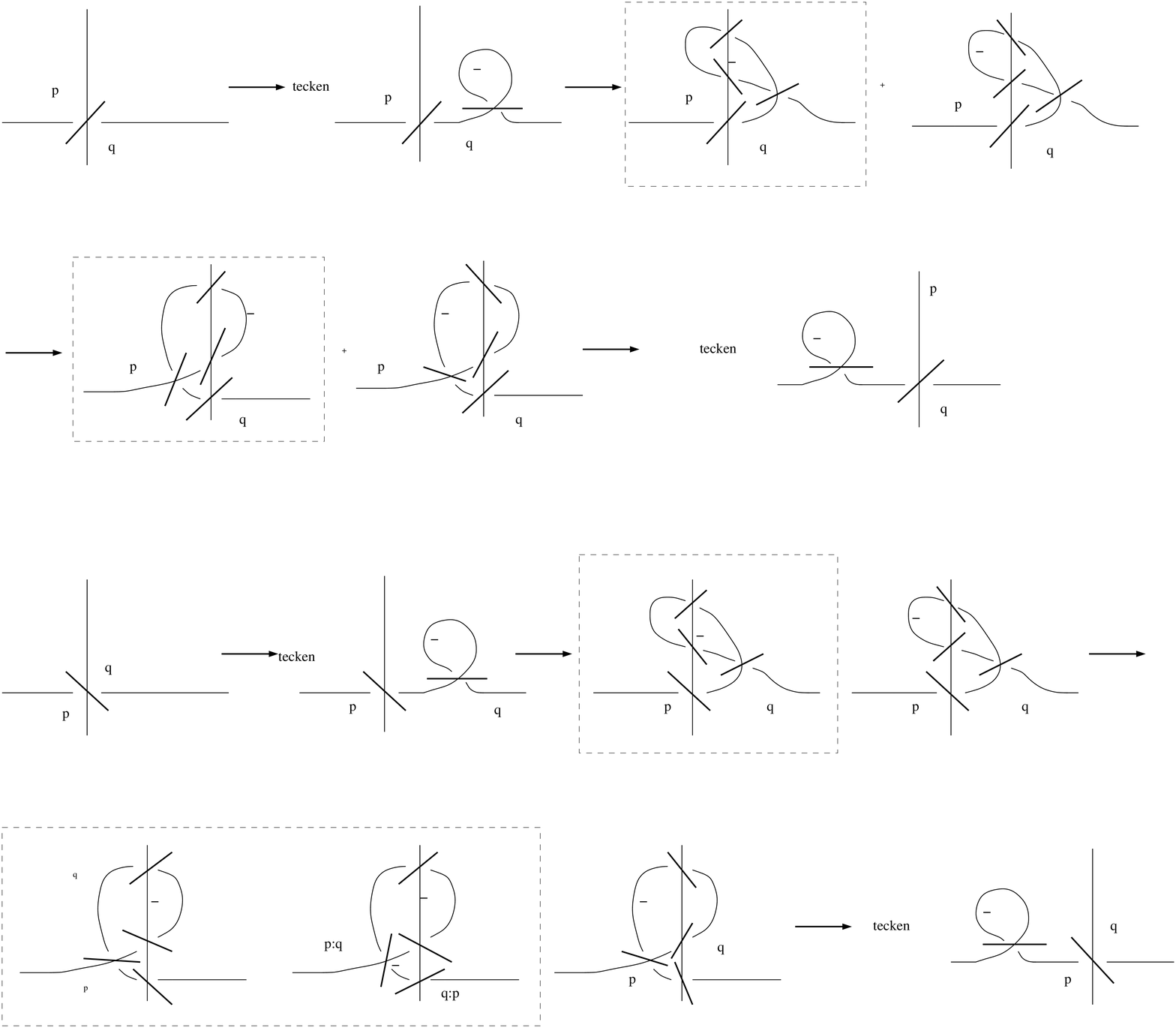}
\caption{Move 6}
\label{move6}
\end{center}
\end{figure}

{\bf Negative twist}\qua
The right side is just the appearance of a negative twist. 
The induced map is given by Figure \ref{chain1}, top.
The left hand side is described in Figure \ref{move6}.

\begin{figure}[ht!]\small
\begin{center}
\psfrag{a}{\rlap{\kern-2pt$a$}}
\psfrag{b}{\rlap{\kern-2pt$b$}}
\psfrag{c}{\rlap{\kern-2pt$c$}}
\psfrag{d}{\rlap{\kern-2pt$d$}}
\includegraphics[width = 12.5 cm, height = 3 cm]{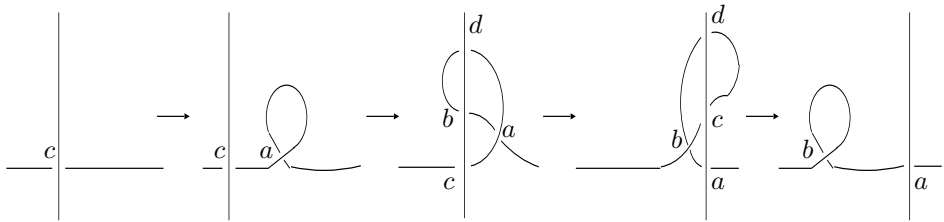}
\caption{Enumeration of the crossings of the left side of Move 6}
\label{aname6}
\end{center}
\end{figure}

The two end results seem to coincide, but we still have to analyze 
the possible sign difference. That is, what happens to the $E(L)$ tensor 
factor. Again, on the right hand side this is given by 
Figure \ref{chain1}. As for the left hand side, pick an enumeration 
of the crossings as in Figure \ref{aname6} and then consider Figure \ref{move6} 
again. The states in dashed boxes eventually map to zero. Therefore 
let us disregard them. Then the $E(L)$--factor
transforms as below for the two types of states (negative respectively 
positive marker at the crossing).
\begin{equation*}
[xc] \mapsto [xca] \mapsto [xcab] \mapsto [xcab] \mapsto [xab]
\end{equation*}
\begin{equation*}
[x] \mapsto [xa] \mapsto [xab] \mapsto [xbc] \mapsto [xb]
\end{equation*}
We see that the signs are the same for both sides.

{\bf Positive twist}\qua In the positive twist case, the two induced
chain maps do not coincide. However, their sum evaluated on a state 
is an element of the contractible subcomplex corresponding to the 
twist in the target diagram. The projection on the non-contractible 
part gives an isomorphism on homology. Tracing the signs as above shows 
that there is a difference in signs here. We omit details and note 
only that the $\Omega_3$--move involved here is the ``reflected'' one 
$\overline{\Omega}_3$ from Section \ref{yetanothermove}.

These remarks apply regardless of whether the horizontal strand is
behind or in front of the vertical one.

\subsubsection*{Movie move 7} 
{\bf Reduction to one version of the move}\qua
It is sufficient to consider a single version of this move, 
for the following reason. Let $P$ be the three-dimensional 
configuration of four planes 
making up the left side of any quadruple point move. (A small 3--ball 
in which four sheets of the surface diagram form a tetrahedron.)
Enumerate the four sheets in the order of increasing height 
with respect to the projection from four to three dimensions. 
Let $P_0$ be the standard configuration which has the $xy$--, 
$xz$-- and $yz$--planes as its first, second and third sheets, and with 
the fourth plane having normal $(1,1,1)$. 

The first, second and third sheets of $P$ can be isotoped to coincide 
with the first, second and third sheets of $P_0$ respectively. 
With this done, the fourth sheet of $P$ can be made to sit as 
one of the eight planes $(\pm 1,\pm 1,\pm 1)$. The moves which 
correspond to isotopies of the surface diagram are the moves 8--15 
and the interchanges of distant critical points. The invariance 
under these moves is proved in other subsections (independently of
the invariance under this move). Thus, this isotopy does not change 
the induced homomorphism (up to sign).

The quadruple point move corresponding 
to this new plane configuration can be replaced with a sequence of three 
movie moves consisting of one move of type $5$ (adding two triple 
points in the diagram), one quadruple point move in a neighbouring 
octant, and another (inverse) move of type $5$ (subtracting two triple points). 
(This is a four-dimensional counterpart of Figure \ref{turaev}, where
one type of $\Omega_3$--move was replaced by a sequence containing
one $\Omega_2$--move, the other $\Omega_3$--move and an inverse 
$\Omega_2$--move.) This result, and pictures explaining it in detail, 
can be found in \cite{CJKLS}, page $11$.
From the invariance under the moves of types $3$ and $5$ we may therefore 
assume that the move takes place in any octant we like, in particular with 
the standard configuration of sheets. 

Finally, the result can be isotoped to the right hand side of the 
original move. Note that orientations never enter into the calculations. 
It follows that invariance under any version of move $7$ is implied 
by the invariance under a single version. We choose 
one with its crossings as in Figure \ref{enumrose7} below, and label the 
crossings as indicated there. The left movie is the upper
one in this figure. 

\begin{figure}[ht!]\small
\begin{center}
\psfrag{a}{\smash{\rlap{\kern-2pt\raise-1pt\hbox{$a$}}}}
\psfrag{b}{\smash{\rlap{\kern-2pt\raise-1pt\hbox{$b$}}}}
\psfrag{c}{\smash{\rlap{\kern-2pt\raise-1pt\hbox{$c$}}}}
\psfrag{d}{\smash{\rlap{\kern-2pt\raise-1pt\hbox{$d$}}}}
\psfrag{e}{\smash{\rlap{\kern-2pt\raise-1pt\hbox{$e$}}}}
\psfrag{f}{\smash{\rlap{\kern-2pt\raise-1pt\hbox{$f$}}}}
\psfrag{D}{\rlap{\kern-2pt$D$}}
\psfrag{D'}{\rlap{\kern-2pt$D'$}}
\includegraphics[width = 10 cm, height = 6 cm]{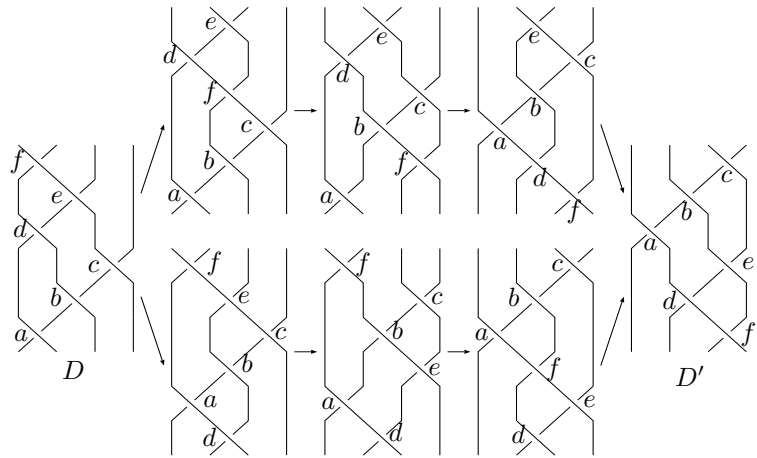}
\caption{Enumeration of the crossings in Move 7}
\label{enumrose7}
\end{center}
\end{figure}

{\bf Computing the maps}\qua
Consider the target diagram $D'$ in Figure \ref{enumrose7}. The bottom right 
triangle, formed by the vertices d,e,f is the target of a third Reidemeister 
move and hence defines a splitting of the chain complex $C(D')$ 
as explained in Section \ref{rmm}, Figures \ref{psi3} and \ref{contr3}. 
In the discussion that follows we will tacitly disregard all states 
in this subcomplex in the discussion that follows. Indeed, projecting 
out this subcomplex induces an isomorphism on homology and it is enough 
to consider the composition of right and left maps with the projection. 

Similarily, the complex of the source diagram $D$ splits according to the upper 
left triangle formed by the vertices d,e,f, where the first $\Omega_3$--move 
in the left movie takes place. We need to consider only the generators 
of the non-contractible factor of this splitting, since the inclusion 
of this factor induces an isomorphism on homology.

\begin{figure}[ht!]\small
\begin{center}
\psfrag{a}{\smash{\rlap{\kern-1.5pt\raise-1pt\hbox{$a$}}}}
\psfrag{b}{\smash{\rlap{\kern-1.5pt\raise-1pt\hbox{$b$}}}}
\psfrag{c}{\smash{\rlap{\kern-1.5pt\raise-1pt\hbox{$c$}}}}
\psfrag{d}{\smash{\rlap{\kern-1.5pt\raise-1pt\hbox{$d$}}}}
\psfrag{e}{\smash{\rlap{\kern-1.5pt\raise-1pt\hbox{$e$}}}}
\psfrag{f}{\smash{\rlap{\kern-1.5pt\raise-1pt\hbox{$f$}}}}
\psfrag{L}{$L$}
\psfrag{R}{$R$}
\includegraphics[width = 3 cm, height = 3 cm]{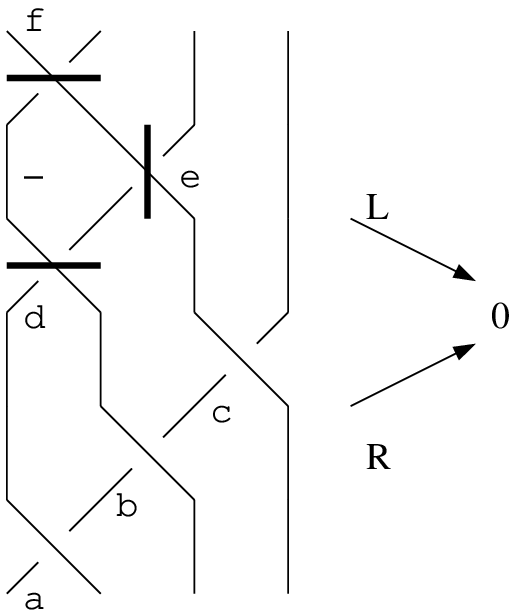}
\caption{Move 7}
\label{a7tilde}
\end{center}
\end{figure}

It is easily checked that a state as in Figure \ref{a7tilde} maps to zero under the 
left hand as well as under the right hand map, regardless of the markers at $a,b,c$. 

It follows that we need only look at the states 
(let us call them $A$--states) which have positive 
markers at $e,f$ and negative at $d$, and all states which 
have a negative marker at $f$ (let us call these $B$--states). 

\begin{figure}[ht!]\small
\begin{center}
\includegraphics[width = 8 cm, height = 3 cm]{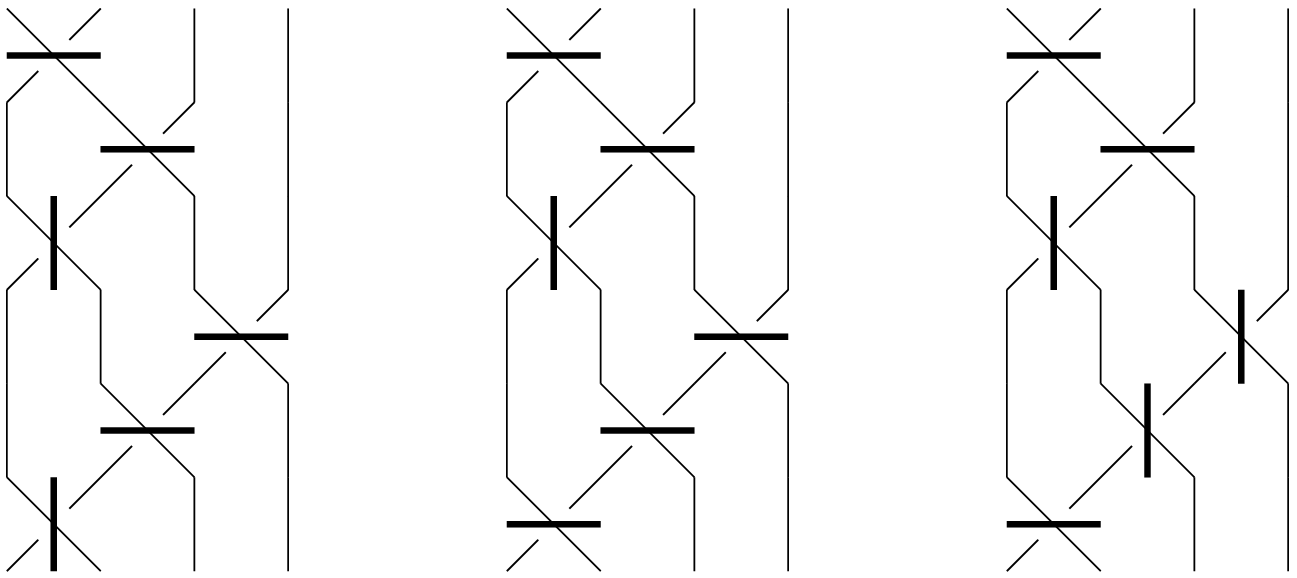}
\caption{Move 7}
\label{a7zero}
\end{center}
\end{figure}

It is also easy to verify that $A$--states with local
configurations of markers as in Figure \ref{a7zero}
all map to zero under both sides of the move.

\begin{figure}[ht!]\small
\begin{center}
\psfrag{+}{$+$}
\psfrag{-}{$-$}
\psfrag{L}{$L$}
\psfrag{R}{$R$}
\includegraphics[width = 7 cm, height = 3 cm]{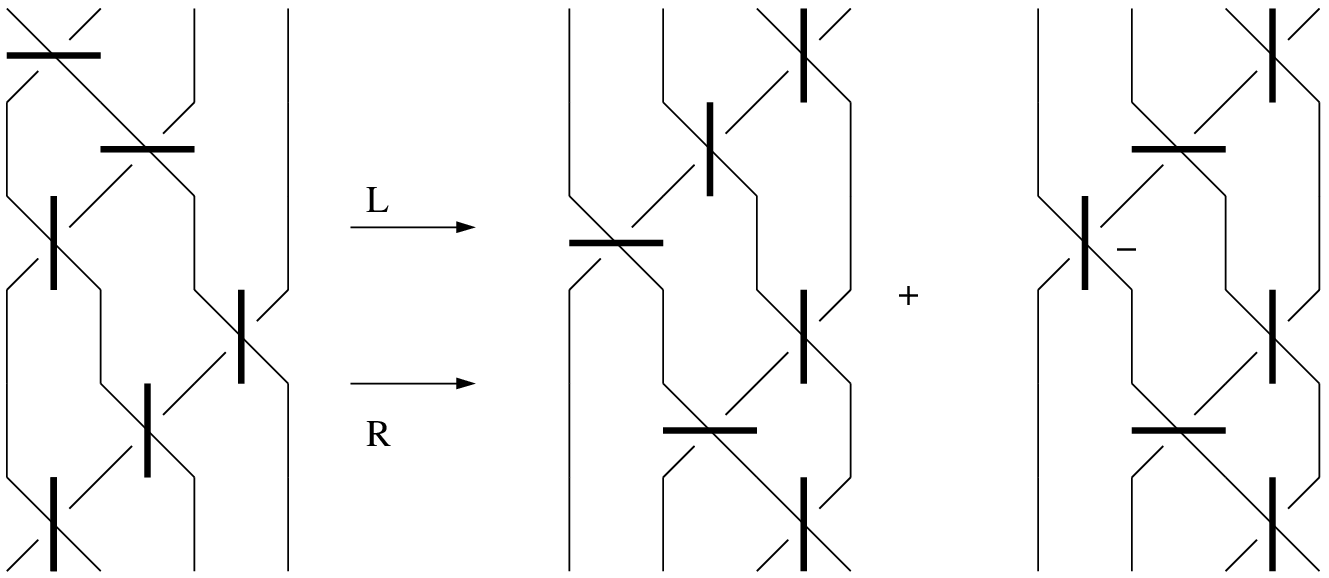}
\caption{Move 7}
\label{a7nnn}
\end{center}
\end{figure}

$A$--states with negative markers at the three remaining crossings
behave in the same way under both sides of the move (Figure \ref{a7nnn}).

We omit here most of the signs of the 
resolutions, but a closer examination \cite{Ja1} reveals that the same
modifications of the signed resolutions take place on both sides, 
so that the final results indeed coincide.

It is tedious but straightforward to verify the results in Figures
\ref{a7pppnp1} and \ref{a7pppnp2} below. Most of the information 
about the signed resolutions is again left out in these figures. 
Consider the expressions in the right column. They are written 
modulo the contractible subcomplex of the bottom right triangle.
The resolutions in the left and right hand images 
are the same (that is, planar isotopic) for two states placed above each other. 
Indeed, they even have the same signed resolutions,
as can be seen by going through the ways signs change under
the two sides of the move.

\begin{figure}[ht!]\small
\begin{center}
\psfrag{+}{$+$}
\psfrag{-}{$-$}
\psfrag{L}{$L$}
\psfrag{R}{$R$}
\includegraphics[width = 11 cm, height = 11 cm]{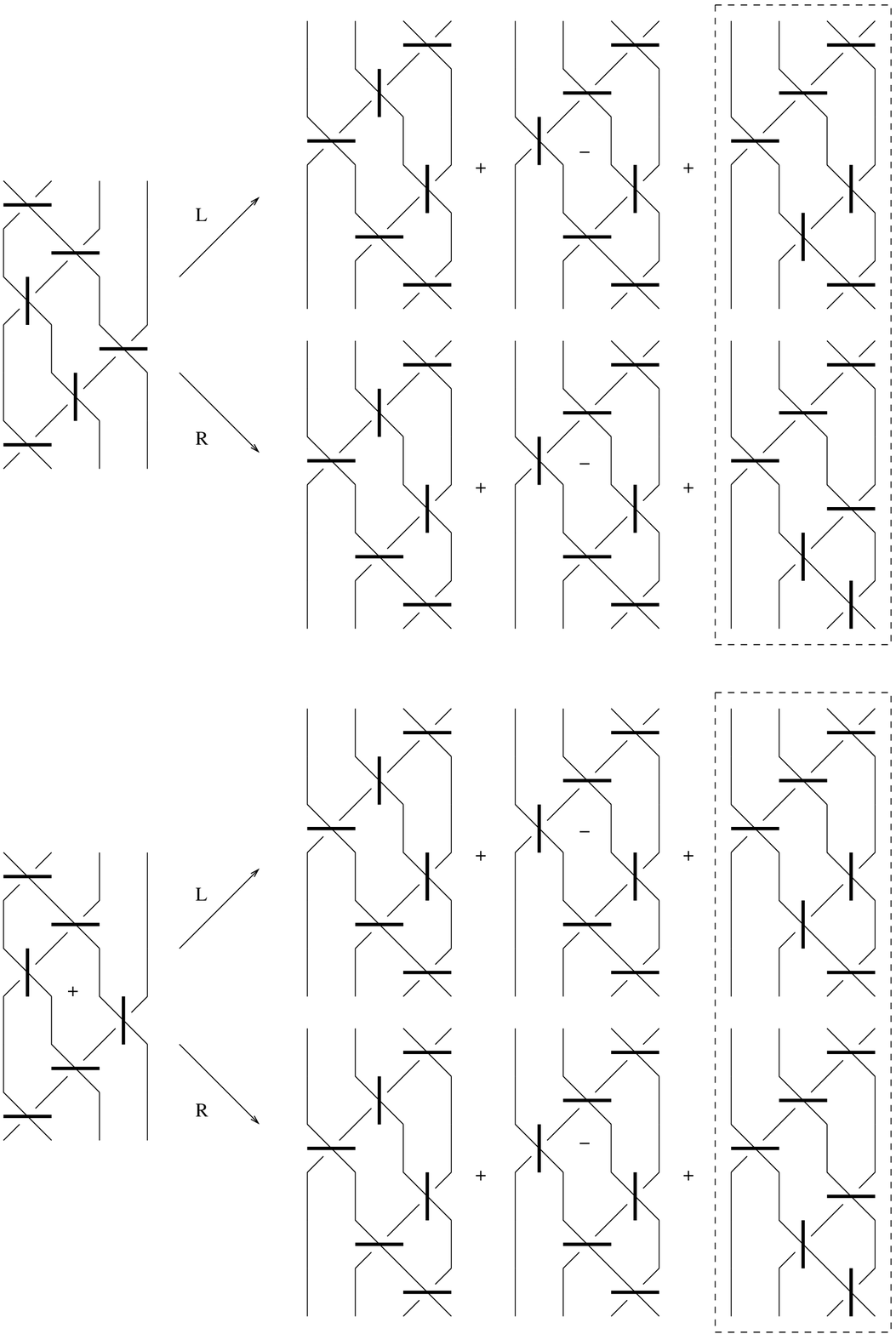}
\caption{Move 7}
\label{a7pppnp1}
\end{center}
\end{figure}

\begin{figure}[ht!]\small
\begin{center}
\psfrag{+}{$+$}
\psfrag{-}{$-$}
\psfrag{L}{$L$}
\psfrag{R}{$R$}
\includegraphics[width = 11 cm, height = 5 cm]{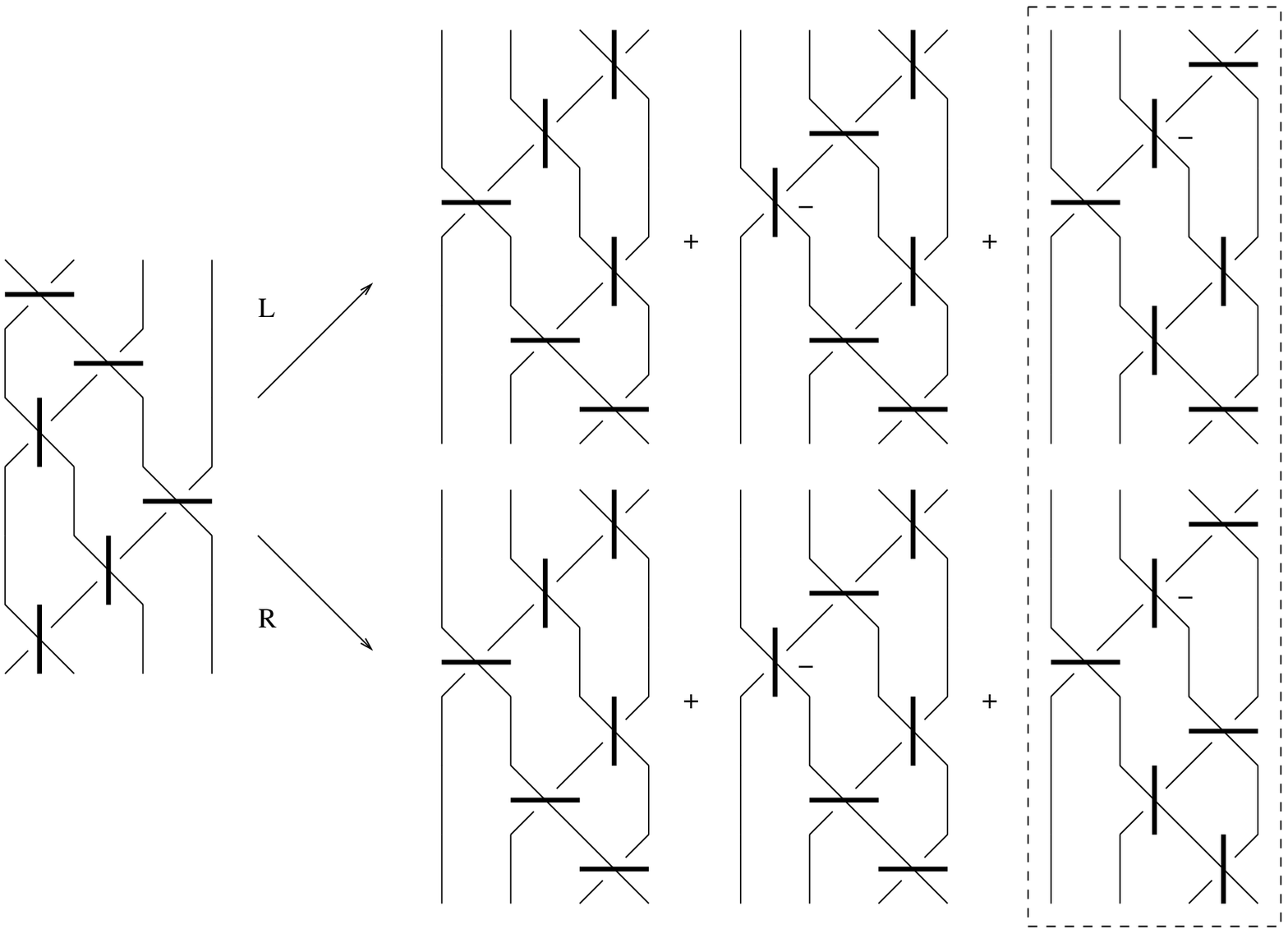}
\includegraphics[width = 13 cm, height = 5 cm]{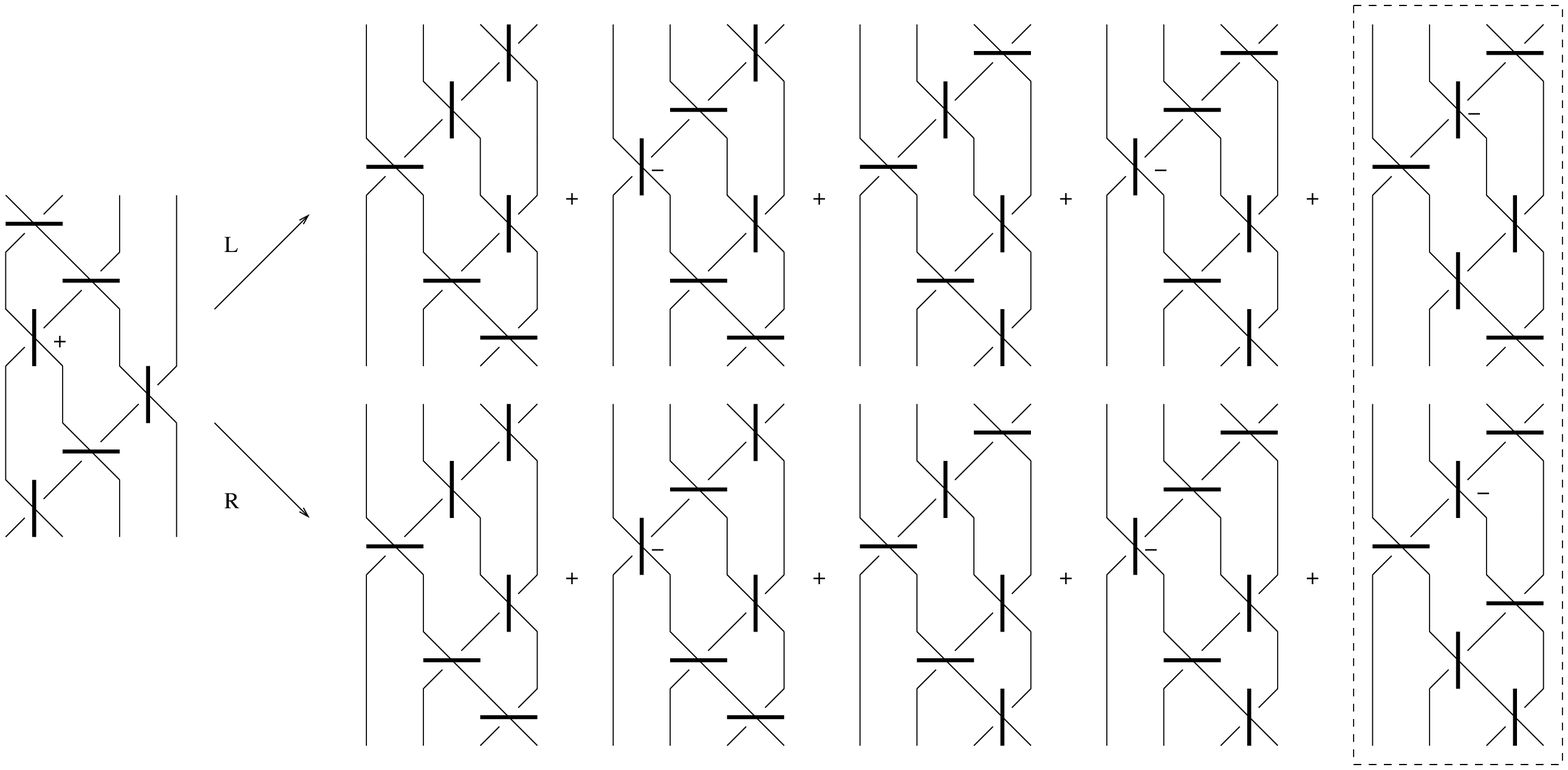}
\caption{Move 7}
\label{a7pppnp2}
\end{center}
\end{figure}

\begin{figure}[ht!]\small
\begin{center}
\psfrag{+}{$+$}
\psfrag{-}{$-$}
\psfrag{p}{$p$}
\psfrag{q}{$q$}
\psfrag{r}{$r$}
\psfrag{tensorxd}{$\otimes [xd])$}
\psfrag{tensorxde}{$\otimes [xde]$}
\psfrag{tensorxdf}{$\otimes [xdf]$}
\psfrag{tensorxdt}{$\otimes [xdt]$}
\psfrag{SUM}{$\sum_{t} $}
\psfrag{d}{$d($}
\psfrag{=}{$=$}
\includegraphics[width = 13 cm, height = 3 cm]{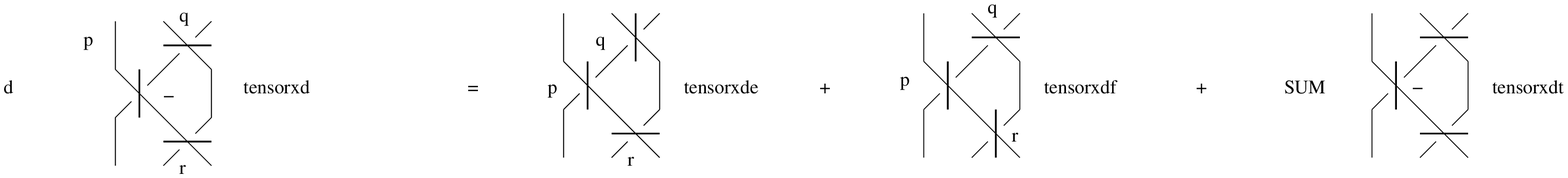}
\caption{Rewriting dashed box differences}
\label{newintegral}
\end{center}
\end{figure}

The next thing to note is that the difference of two states 
enclosed in a dashed box is really in the contractible subcomplex. 
Such states differ only at the lower right triangle 
and as mentioned above they have the same signed resolutions. That is, 
they look like the first two terms on the right hand side of the
equation in Figure \ref{newintegral}. At first glance, it seems 
as though we should get the difference rather than the sum 
of these terms. However, a closer analysis reveals the correct 
sign hidden in the $E(L)$--factor. Since all terms except for these two 
are in $C_{contr}$, the statement follows.

Finally, consider again $A$--states with markers (pos, pos, neg)
respectively (neg, pos, neg) at crossings (a, b, c), but this time 
with a small negative component (the circle $bced$) instead of a
positive one. On these states both sides are zero as a short 
computation shows. This finishes the proof for $A$--states,
except for the verification that the $E(L)$--factor really 
behaves as claimed. Using the enumeration of vertices 
above, this is a straightforward check.

As regards the $B$--states, i.e.\ states with a negative marker at $f$,
these are more well-behaved. Indeed, the left and right images
are the same for each of these states, so the maps are the same on the
chain level, already. We omit details.

\subsubsection*{Movie move 8}
\begin{figure}[ht!]\small
\begin{center}
\psfrag{rp}{$)$}
\psfrag{L}{$L$}
\psfrag{R}{$R$}
\psfrag{p}{$p$}
\psfrag{tecken}{$(-1)^i$(}
\psfrag{-}{$-$}
\psfrag{+}{$+$}
\includegraphics[width = 10 cm, height = 5 cm]{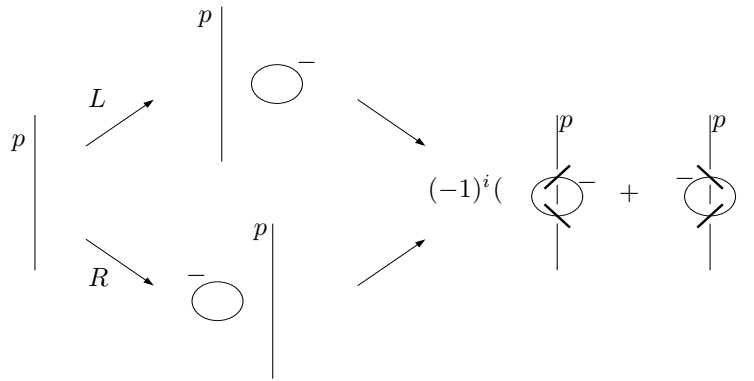}
\caption{Move 8 in downward time}
\label{a8down}
\end{center}
\end{figure}

In  Figure \ref{a8down}, we give the computation for the 
move no. $8$ with time flowing downwards. A new-born circle 
is always marked with a minus sign, and the $\Omega_2$--move
that follows eliminates the difference between the two sides.

\begin{figure}[ht!]\small
\begin{center}
\psfrag{tecken}{$(-1)^i$}
\psfrag{teckenlp}{$(-1)^i$(}
\psfrag{tecken+1}{$(-1)^{i+1}$}
\psfrag{rp}{$)$}
\psfrag{-}{$-$}
\psfrag{+}{$+$}
\psfrag{R}{$L$}
\psfrag{L}{$R$}
\psfrag{0}{0}
\includegraphics[width = 12 cm, height = 15 cm]{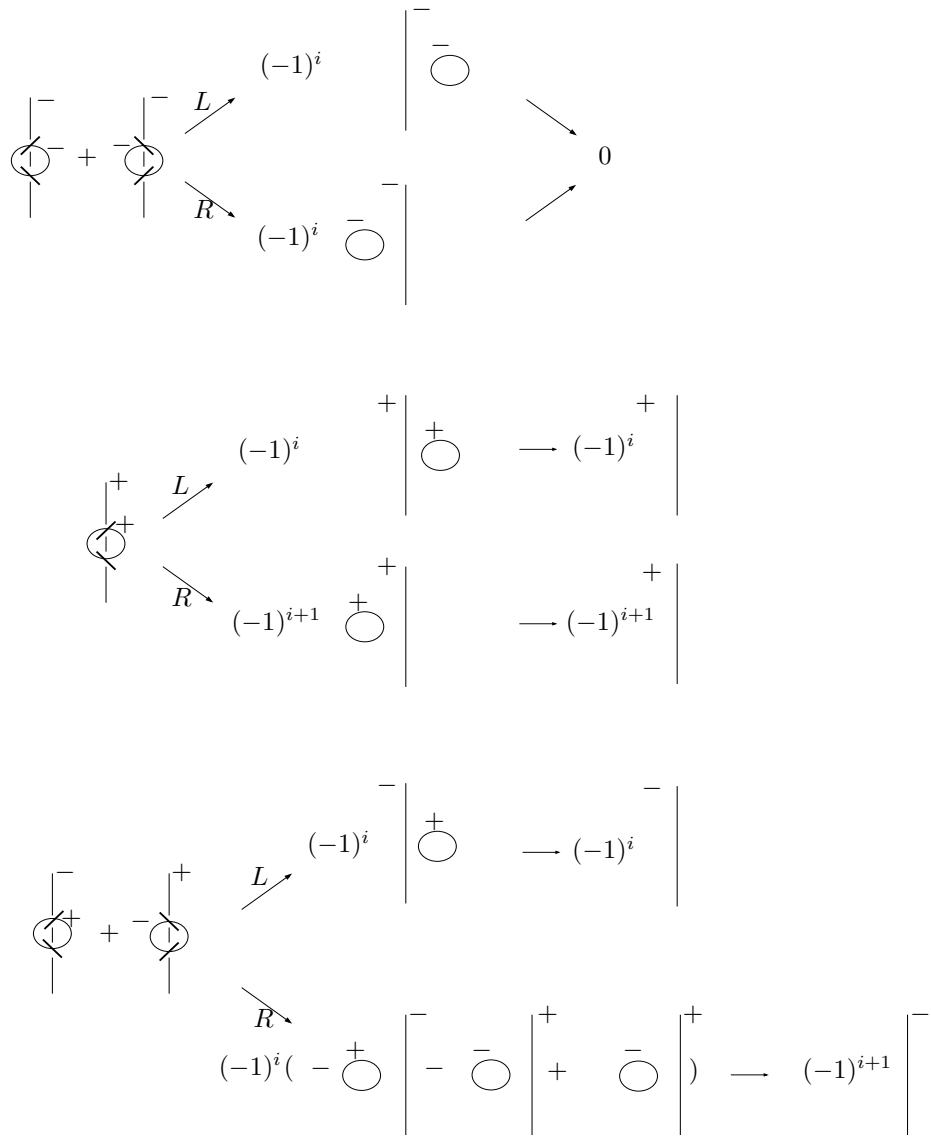}
\caption{Move 8 in upward time}
\label{a8up}
\end{center}
\end{figure}

In upward time, we can start by splitting the chain complex into
$C \oplus C_{contr}$ according to the $\Omega_2$--move about to be 
performed on the left hand side. It is enough to check what happens 
to the generators of $C$ (see Figure \ref{psi2}). The calculation 
is included in Figure \ref{a8up}. 

No extra signs appear from the $E(L)$--factor, neither in upward nor in
downward time. The version where the vertical strand is above the 
circle is almost identical.

\subsubsection*{Movie move 9}
\begin{figure}[ht!]\small
\begin{center}
\psfrag{+}{$+$}
\psfrag{-}{$-$}
\psfrag{p}{$p$}
\includegraphics[width = 10 cm, height = 5 cm]{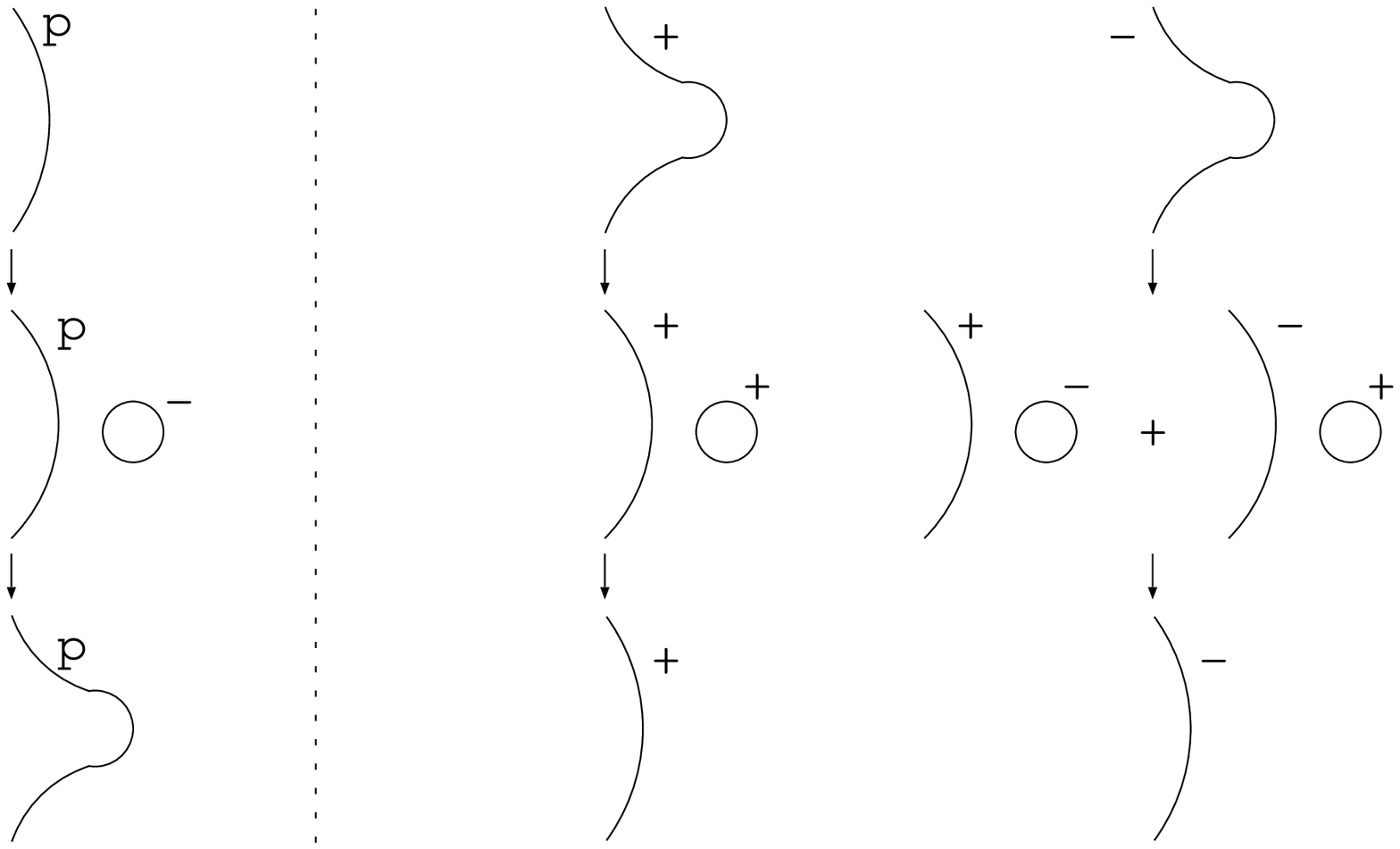}
\caption{Move 9}
\label{maps9}
\end{center}
\end{figure}

The computation here is very trivial on both sides of the move. 
Only maps associated to Morse modifications are used.
We immediately see that the maps coincide up to sign.
No markers are involved, so obviously the effect on the 
$E(L)$--factor is trivial. See Figure \ref{maps9}.

\subsubsection*{Movie move 10} 
\begin{figure}[ht!]\small
\begin{center}
\psfrag{+}{$+$}
\psfrag{-}{$-$}
\psfrag{p}{$p$}
\psfrag{p:q}{$p{:}q$}
\psfrag{q:p}{$q{:}p$}
\psfrag{q:r}{$q{:}r$}
\psfrag{r:q}{$r{:}q$}
\psfrag{q}{$q$}
\psfrag{r}{$r$}
\psfrag{tecken}{$(-1)^i$(}
\psfrag{rp}{$)$}
\includegraphics[width = 13 cm, height = 7 cm]{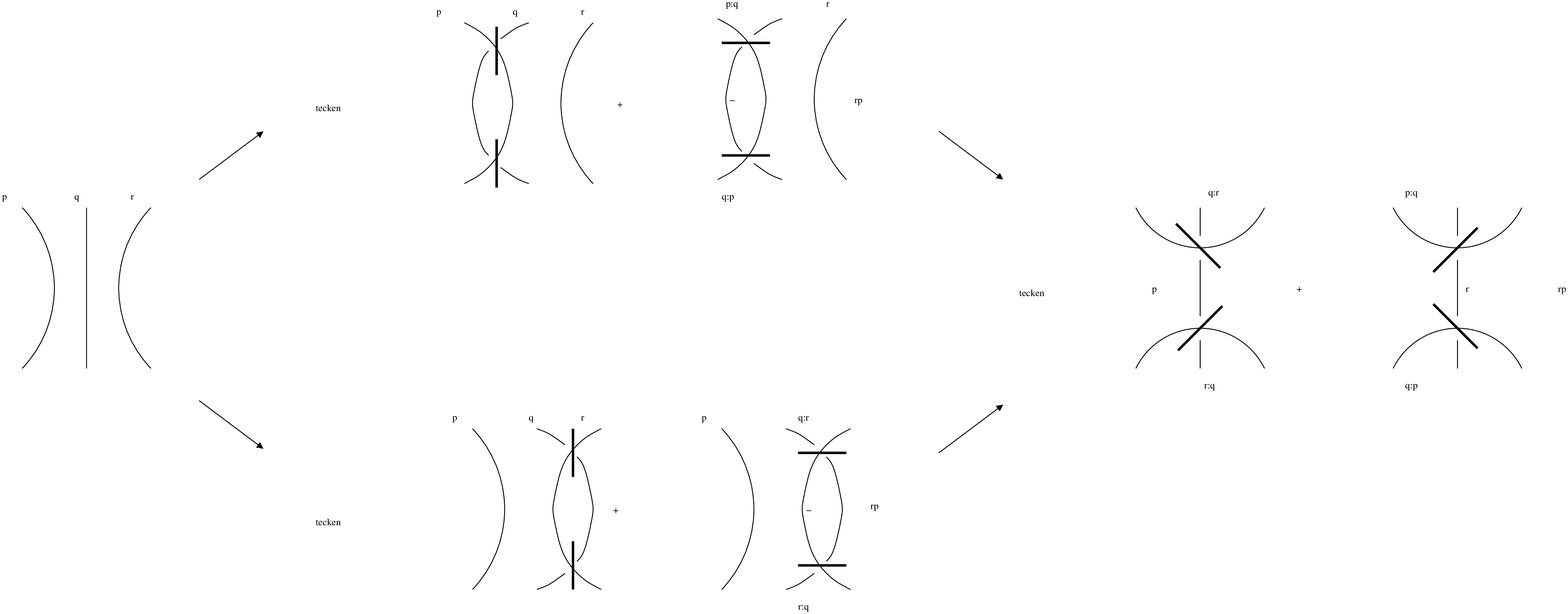}
\caption{Move 10 in downward time}
\label{a10down}
\end{center}
\end{figure}

In downward time the calculation in 
Figure \ref{a10down} yields the result. 

\begin{figure}[ht!]\small
\begin{center}
\psfrag{+}{$+$}
\psfrag{-}{$-$}
\psfrag{p:q}{$p{:}q$}
\psfrag{q:p}{$q{:}p$}
\psfrag{q:r}{$q{:}r$}
\psfrag{r:q}{$r{:}q$}
\psfrag{p}{$p$}
\psfrag{q}{$q$}
\psfrag{r}{$r$}
\psfrag{tecken}{$(-1)^i$}
\psfrag{tecken+1}{$(-1)^{i+1}$}
\includegraphics[width = 8 cm, height = 12 cm]{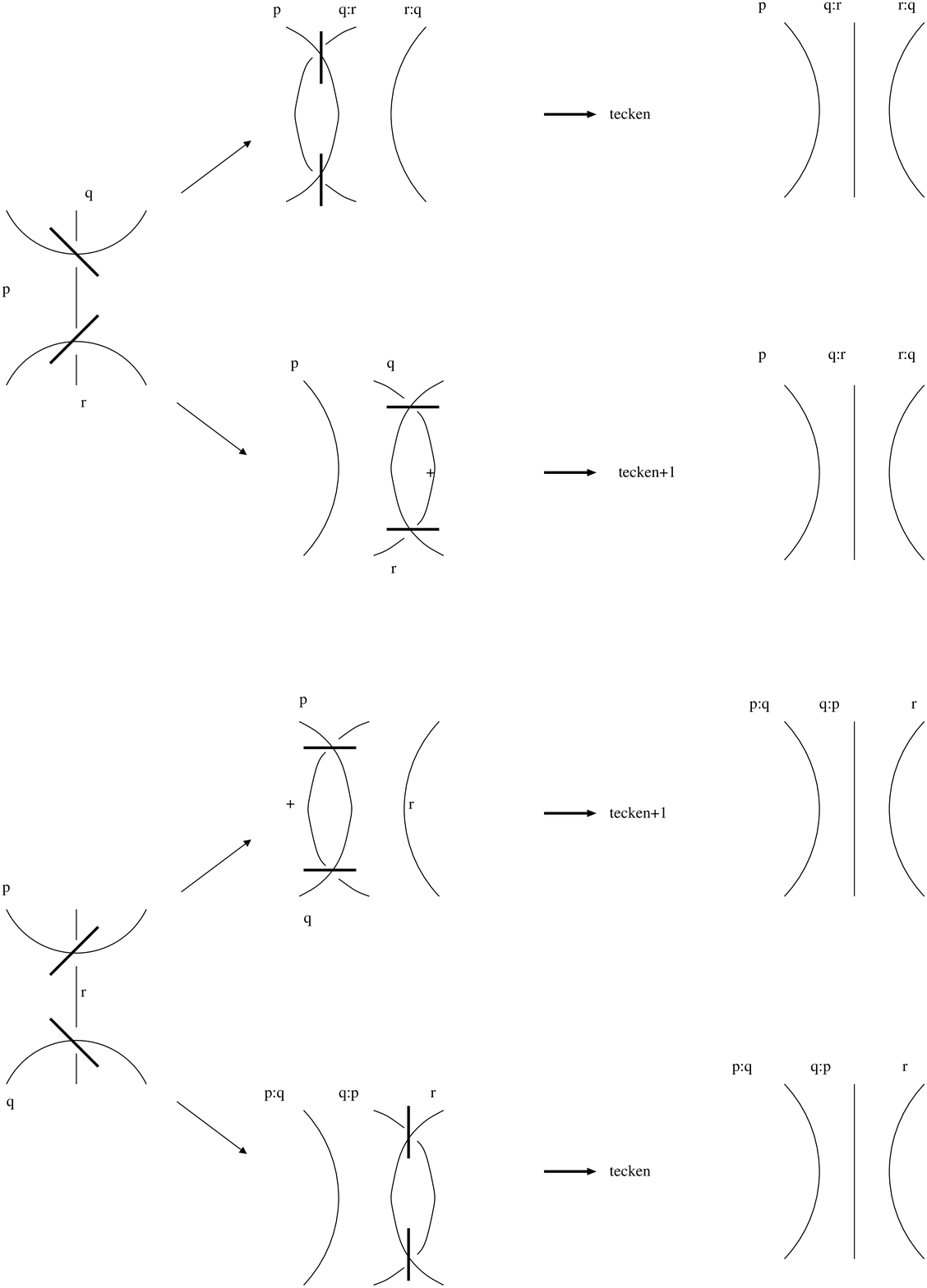}
\caption{Move 10 in upward time: on the second level terms 
in $C_{contr}$ are not drawn.}
\label{a10up}
\end{center}
\end{figure}

And in upward time, the one in Figure \ref{a10up}, in which one
should note that the second level is written modulo $C_{contr}$. The states 
not shown all map to zero.

The action on the $E(L)$--factor is the same on both sides, and in both
time directions. If the middle strand is behind the others,
everything works similarily. 

\subsubsection*{Movie moves 11 -- 15}
The calculations for the rest of the moves are similar to the above
and straightforward (albeit sometimes tedious) given Proposition \ref{prop}.
We omit them, and refer the interested reader to \cite{Ja1} for details.
There the reader can also find more details for the previous
cases, if (s)he so wishes. \qed

\begin{lemma}The invariant is non-trivial. 
\end{lemma}
\begin{proof}
Let $K$ be a knot with non-trivial Jones polynomial. (No non-trivial
knots with trivial Jones polynomial are known at present.)
For any knot $K$ in $\mathbb{R}³$, it is well-known that the connected 
sum of $K$ and its mirror image $\bar{K}$ is ``slice'', i.e.\ bounds 
an embedded disc in $\mathbb{R}^3 \times I$. Put one such disc at each end of
$\mathbb{R}^3 \times I$ and connect the two discs with a vertical tube.
This is a knotted cylinder and induces some map on homology.
In general, this map has a big kernel, since it factors 
through the homology of the unknot in the space between the two discs.
Therefore it is different from the identity cylinder on $K \# \bar{K}$
if the latter has non-trivial homology groups. 
This proves non-triviality, and concludes the whole proof. 
\end{proof}


\section{Lefschetz polynomials of link endocobordisms}
\label{lple}
The calculation of the induced homomorphism on homology induced 
by a link cobordism can be rather involved. However it 
is relatively easy to compute it on the level of chain 
groups, because of their explicit, geometric definition. 
It would therefore be interesting to find (possibly weaker)
derived invariants of link cobordisms computable from the 
homomorphism induced on the chain group level.

In this section we make some remarks on special link cobordisms, namely
those which have identical source and target $L$ and the cobordism $\Sigma$ 
a surface with zero Euler characteristic. We could call such objects 
$(\Sigma, L)$ ``links with endomorphism'' or ``link endocobordisms''. 

By Theorem \ref{mkc}, to every link endocobordism $(\Sigma, L)$ and diagram $D$ of $L$ 
there is a well-defined (up to sign) endomorphism $\phi_* = \phi_{\Sigma *}$ 
on $\mathcal{H}(D)$. On $\mathcal{H}^{i,j}(D)$ we denote this map by $\phi_*^{i,j}$. Note
that the  bigrading is $(0,0)$. For fixed $j$ this is an endomorphism 
of a chain complex, so it has a well-defined Lefschetz number
$L_{j}(\Sigma)$ (also up to sign) which is the alternating sum of traces: 
\begin{displaymath}
L_{j}(\Sigma) = \sum_i (-1)^i tr(\phi_*^{i,j}:\mathcal{H}^{i,j}(D) \otimes
\mathbb{Q}\rightarrow \mathcal{H}^{i,j}(D) \otimes \mathbb{Q})
\end{displaymath}
$L_j(\Sigma)$ is obviously invariant (up to sign) under ambient
isotopy, as is the following alternating sum.

\begin{definition}
Summing over $j$ with coefficients $q^{j}$ gives us {\em the Lefschetz
polynomial} of the link endocobordism:
\begin{displaymath}
L(\Sigma) = \sum_{j} L_{j}(\Sigma)q^j
\end{displaymath}
\end{definition}

Isotopies of $L$ act in a natural way on endocobordisms $(\Sigma,L)$ 
with diagram $D$. Namely, to each isotopy 
of $L$ corresponds a sequence of Reidemeister moves on $D$, 
that is, a  movie from $D$ to another diagram $D'$. 
Conjugating $\Sigma$ by the corresponding link cobordism, we get a new pair 
$(\Sigma',L')$ with diagram $D'$.
Since the Lefschetz polynomial is built from traces it is invariant 
under conjugation. In fact, even each $L_j$ is invariant.

There is the well-known theorem that $L_j$ can be computed 
directly on the chain level, without passing to homology, so that 
\begin{displaymath}
L_{j}(\Sigma) = \sum_i (-1)^i tr(\phi^{i,j}:C^{i,j}(D) \otimes
\mathbb{Q} \rightarrow C^{i,j}(D) \otimes \mathbb{Q}).
\end{displaymath}
where $\phi^{i,j} = \phi^{i,j}_\Sigma$ is the chain map induced by $\Sigma$.

\begin{remark}
It is easy to see that the map $\phi^{i,-7}_\Sigma$ of clockwise rotation 
in Section \ref{cekc}, only has the following 
non-zero matrices in the canonical bases of $C^{i,-7}(D)$.  
\[
\phi^{-4,-7}_{\Sigma} = 
\left(
\begin{array}{ccccc}
0&0&0&1&0\\
1&0&0&0&0\\
0&1&0&0&0\\
0&0&1&0&0\\
0&0&0&0&1\\
\end{array}
\right)
\]
\[
\phi^{-3,-7}_{\Sigma} =
\left(
\begin{array}{cccccccc}
0&0&0&1&0&0&0&0\\
1&0&0&0&0&0&0&0\\
0&1&0&0&0&0&0&0\\
0&0&1&0&0&0&0&0\\
0&0&0&0&0&0&0&1\\
0&0&0&0&1&0&0&0\\
0&0&0&0&0&1&0&0\\
0&0&0&0&0&0&1&0\\
\end{array}
\right)
\]
From this we immediately see that $L_{-7}(\Sigma) = (\pm)1$.
(This can also be seen from the matrix of $\phi_{{\Sigma} *}$ 
in the remark at the end of Section \ref{cekc}.)
By contrast, the identity map has this Lefschetz number equal to
$-3$. Hence the induced maps cannot be the same.
\end{remark}

{\bf Acknowledgements}\qua This paper and its longer predecessor
\cite{Ja1} were written while I was a graduate student at Uppsala
University. The final version was prepared at INdAM, Rome. I thank the
referee for several useful comments.

\bibliographystyle{gtart}

\Addresses\recd
\end{document}